# MULTIVARIATE ARCHIMEDEAN COPULAS, $D$-MONOTONE FUNCTIONS AND $\ell_1$-NORM SYMMETRIC DISTRIBUTIONS

BY ALEXANDER J. MCNEIL AND JOHANNA NEŠLEHOVÁ

*Maxwell Institute, Edinburgh and ETH Zurich*

It is shown that a necessary and sufficient condition for an Archimedean copula generator to generate a $d$-dimensional copula is that the generator is a $d$-monotone function. The class of $d$-dimensional Archimedean copulas is shown to coincide with the class of survival copulas of $d$-dimensional $\ell_1$-norm symmetric distributions that place no point mass at the origin. The $d$-monotone Archimedean copula generators may be characterized using a little-known integral transform of Williamson [*Duke Math. J.* **23** (1956) 189–207] in an analogous manner to the well-known Bernstein–Widder characterization of completely monotone generators in terms of the Laplace transform. These insights allow the construction of new Archimedean copula families and provide a general solution to the problem of sampling multivariate Archimedean copulas. They also yield useful expressions for the $d$-dimensional Kendall function and Kendall's rank correlation coefficients and facilitate the derivation of results on the existence of densities and the description of singular components for Archimedean copulas. The existence of a sharp lower bound for Archimedean copulas with respect to the positive lower orthant dependence ordering is shown.

**1. Introduction.** Archimedean copulas are an important class of multivariate dependence models which first appeared in the context of probabilistic metric spaces (see, e.g., [37]) and which enjoy considerable popularity in a number of practical applications. Following [22], any Archimedean copula has the simple algebraic form

(1.1) $$C(u_1, \ldots, u_d) = \psi(\psi^{-1}(u_1) + \cdots + \psi^{-1}(u_d)),$$
$$(u_1, \ldots, u_d) \in [0,1]^d,$$

where $\psi$ is a specific function known as a generator of $C$. As a consequence, many dependence properties of such copulas are relatively easy to establish because they reduce to analytical properties of the generator $\psi$; see, for example, [11, 12, 15, 19, 26, 27]. Numerous statistically tractable parametric families of Archimedean copulas have been constructed and used in a variety of practical applications, including multivariate survival analysis [5, 30], actuarial loss modeling [9, 21] and quantitative finance [4, 36].

While any $d$-dimensional Archimedean copula is necessarily of the form (1.1), the converse of this statement is not true. Somewhat surprisingly, the conditions under which a generator $\psi$ defines a $d$-dimensional copula by means of (1.1) have not yet been fully clarified except in two cases. Schweizer and Sklar [37] show that a generator $\psi$ induces a bivariate copula if and only if it is convex, whereas Kimberling [20] proves that $\psi$ defines an Archimedean copula in any dimension if and only if it is a completely monotone function, or equivalently, a Laplace transform of a nonnegative random variable. Kimberling's condition is, however, not necessary for a given dimension $d \geq 3$ and leads to limited dependence characteristics. Some authors [14, 26, 28] have required that $\psi$ only has derivatives up to order $d$ which alternate in sign, a condition which, though sufficient and considerably weaker, is still not necessary. The present paper fills this gap by showing that the necessary and sufficient condition is that $\psi$ should have an analytical property known as $d$-monotonicity. This allows the possibility of $d$-dimensional Archimedean copulas without densities and reveals the existence of a sharp lower bound on the set of all $d$-dimensional Archimedean copulas with respect to the positive lower orthant dependence ordering as defined in [19].

The understanding of the necessary and sufficient conditions for $\psi$ leads to the insight that Archimedean copulas have a very natural geometric interpretation: they are the copulas that are implicit in the survival functions of $d$-dimensional $\ell_1$-norm symmetric distributions, also known as simplicially contoured distributions. These were first introduced by Fang and Fang [7] and comprise scale mixtures of the uniform distribution on the unit $\ell_1$-norm sphere. While it has been pointed out by Alfred Müller that Archimedean copulas can be constructed using $\ell_1$-norm symmetric distributions, this paper shows that these distributions are, in fact, absolutely central to the study of Archimedean copulas.

Formalization of this link requires the consideration of a little-known integral transform which appears in the work of Williamson [41] and which we call in this paper a Williamson $d$-transform. Essentially this transform plays the same role in the study of $d$-monotone generators that the Laplace transform plays in the study of completely monotone generators and knowledge of this link enables us both to propose a general solution to the problem of generating Archimedean copulas with arbitrary given generators and to construct and analyze a rich variety of new Archimedean copulas.

The paper is organized as follows: Section 2 discusses the necessary and sufficient conditions for a function $\psi$ to generate a $d$-dimensional Archimedean copula. Section 3 establishes the connections between Archimedean copulas, $\ell_1$-norm symmetric distributions and Williamson $d$-transforms. The theory is then utilized in Section 4 to prove several useful properties of Archimedean copulas. This includes establishing conditions for the existence of densities and singular components, extending the notion of the Kendall function (a cornerstone of nonparametric inference for bivariate Archimedean copulas introduced by Genest and Rivest [15]) to higher dimensions, and deriving a lower bound for Archimedean copulas with respect to the positive lower orthant dependence ordering. In Section 5 we comment on implications for stochastic simulation of Archimedean copulas and statistical inference.

**2. Archimedean copulas in a given dimension.** This section gives necessary and sufficient conditions under which a generator $\psi$ induces an Archimedean copula via (1.1). To begin, we need to establish basic notation. Throughout, $\mathbf{x}$ denotes a vector $(x_1, \ldots, x_d)$ in $\mathbb{R}^d$; in particular, $\mathbf{0}$ is the origin. If not otherwise stated, all expressions such as $\mathbf{x} + \mathbf{y}$, $\max(\mathbf{x}, \mathbf{y})$ or $\mathbf{x} \leq \mathbf{y}$ are understood as componentwise operations. Furthermore, $[\mathbf{x}, \mathbf{y}]$ refers to the set $[x_1, y_1] \times \cdots \times [x_d, y_d]$ and $\mathbb{R}_+^d$ abbreviates the positive quadrant $[0, \infty)^d$. Finally, $\|\mathbf{x}\|_1$ denotes the $\ell_1$-norm of $\mathbf{x}$, that is, $\|\mathbf{x}\|_1 = \sum_{i=1}^d |x_i|$, and $x_+$ denotes $\max(x, 0)$.

The symbol $\mathbf{X}$ will be reserved for a random vector on $\mathbb{R}^d$ with distribution function $H$ and survival function $\bar{H}$, defined, respectively, by

$$H(\mathbf{x}) = P(\mathbf{X} \leq \mathbf{x}) \quad \text{and} \quad \bar{H}(\mathbf{x}) = P(\mathbf{X} > \mathbf{x})$$

for any $\mathbf{x} \in \mathbb{R}^d$. Note that $\bar{H}$ may be uniquely retrieved from $H$ by means of the Sylvester–Poincaré sieve formula and conversely. If the support of $H$ is a subset of $(0, \infty)^d$, $\bar{H}$ is uniquely given by its values in $\mathbb{R}_+^d$ because obviously $\bar{H}(\mathbf{x}) = \bar{H}(\max(\mathbf{x}, \mathbf{0}))$, in which case we refrain from specifying the values of $\bar{H}$ outside of $\mathbb{R}_+^d$. Note, however, that if the associated probability distribution places mass on $\mathbb{R}^d \setminus (0, \infty)^d$, in particular on the boundary of $\mathbb{R}_+^d$, $\bar{H}$ *cannot* be uniquely recovered from its restriction to $\mathbb{R}_+^d$.

Finally, we will use *difference operators* defined as follows: Let $f$ be an arbitrary $d$-place real function, $\mathbf{x} \in \mathbb{R}^d$ and $\mathbf{h} > \mathbf{0}$. Then the $d$th order difference $\Delta_{\mathbf{h}} f(\mathbf{x})$ is defined as

$$\Delta_{\mathbf{h}} f(\mathbf{x}) = \Delta_{h_d}^d \cdots \Delta_{h_1}^1 f(\mathbf{x}),$$

where $\Delta_{h_i}^i$ denotes the first-order difference operator given by

$$\Delta_{h_i}^i f(\mathbf{x}) = f(x_1, \ldots, x_{i-1}, x_i + h_i, x_{i+1}, \ldots, x_d) \\ - f(x_1, \ldots, x_{i-1}, x_i, x_{i+1}, \ldots, x_d).$$

In the context of distribution functions, $\Delta_{\mathbf{h}}H(\mathbf{x})$ is the volume assigned by $H$ to the interval $(\mathbf{x}, \mathbf{x}+\mathbf{h}]$. This motivates us to define a function $f : A \to \mathbb{R}$, $A \subseteq \mathbb{R}^d$ fulfilling $\Delta_{\mathbf{h}} f(\mathbf{x}) \geq 0$ for any choice of $\mathbf{x}$ and $\mathbf{h}$ so that all vertices of $(\mathbf{x}, \mathbf{x}+\mathbf{h}]$ that lie in $A$ as *quasi-monotone* on $A$. In other accounts quasi-monotone functions are referred to as $d$-increasing, but we prefer the former terminology here to avoid confusion with the notion of $d$-monotonicity of real functions, which will play a key role later on.

Next, we recall few basic results which will be needed in subsequent discussions.

DEFINITION 2.1. A ($d$-dimensional) copula is a function $C : [0,1]^d \to [0,1]$ satisfying

(i) $C(u_1, \ldots, u_d) = 0$ whenever $u_i = 0$ for at least one $i = 1, \ldots, d$.
(ii) $C(u_1, \ldots, u_d) = u_i$ if $u_j = 1$ for all $j = 1, \ldots, d$ and $j \neq i$.
(iii) $C$ is quasi-monotone on $[0,1]^d$.

Perhaps unfortunately, a survey of Archimedean copulas will amount to the investigation of certain probability distributions on $\mathbb{R}^d$ given by their survival rather than their distribution functions. To avoid working with the more cumbersome formula for $H$ in terms of $\bar{H}$, it is convenient to formulate the following simple observation:

LEMMA 1. *A $d$-place function $\bar{H} : \mathbb{R}^d \to [0,1]$ is a survival function of a probability measure on $\mathbb{R}^d$ if and only if*

(i) $\bar{H}(-\infty, \ldots, -\infty) = 1$ *and* $\bar{H}(\mathbf{x}) = 0$ *if* $x_i = \infty$ *for at least one* $i = 1, \ldots, d$.
(ii) $\bar{H}$ *is right-continuous, that is, for all* $\mathbf{x} \in \mathbb{R}^d$ *it holds that*

$$\forall \varepsilon > 0 \; \exists \delta > 0 \; \forall \mathbf{y} \geq \mathbf{x} \quad \|\mathbf{y} - \mathbf{x}\|_1 < \delta \Rightarrow |\bar{H}(\mathbf{y}) - \bar{H}(\mathbf{x})| < \varepsilon.$$

(iii) *The function $G$ given by $G(\mathbf{x}) = \bar{H}(-\mathbf{x})$, $\mathbf{x} \in \mathbb{R}^d$, is quasi-monotone on $\mathbb{R}^d$.*

PROOF. Since the proof is a standard exercise, we only provide a sketch. First, assume $\mathbf{X}$ is a random vector with survival function $\bar{H}$. Then (i) and (ii) are immediate and (iii) is due to the fact that $\bar{H}(-\mathbf{x}) = P(-\mathbf{X} < \mathbf{x})$. Conversely, suppose $\bar{H}$ satisfies conditions (i)–(iii). Then the right-continuous version $G_+$ of $G$ is a distribution function on $\mathbb{R}^d$. Let $\mathbf{X}$ be a random vector with distribution function $G_+$ and observe that $\bar{H}$ is the survival function of $-\mathbf{X}$. □

Similarly, it will prove convenient to restate the original result by Sklar [38] in terms of survival functions.

THEOREM 2.1. *Let $\bar{H}$ be a $d$-dimensional survival function with margins $\bar{F}_i$, $i = 1, \ldots, d$. Then there exists a copula $C$, referred to as the survival copula of $\bar{H}$, such that*

$$\bar{H}(\mathbf{x}) = C(\bar{F}_1(x_1), \ldots, \bar{F}_d(x_d)) \tag{2.1}$$

*for any $\mathbf{x} \in \mathbb{R}^d$. Furthermore, $C$ is uniquely determined on $D = \{\mathbf{u} \in [0,1]^d : \mathbf{u} \in \operatorname{ran}\bar{F}_1 \times \cdots \times \operatorname{ran}\bar{F}_d\}$ where $\operatorname{ran}\bar{F}_i$ denotes the range of $\bar{F}_i$. In addition, for any $\mathbf{u} \in D$,*

$$C(\mathbf{u}) = \bar{H}(\bar{F}_1^{-1}(u_1), \ldots, \bar{F}_d^{-1}(u_d)),$$

*where $\bar{F}_i^{-1}(u_i) = \inf\{x : \bar{F}_i(x) \leq u_i\}$, $i = 1, \ldots, d$. Conversely, given a copula $C$ and univariate survival functions $\bar{F}_i$, $i = 1, \ldots, d$, $\bar{H}$ defined by (2.1) is a $d$-dimensional survival function with marginals $\bar{F}_1, \ldots, \bar{F}_d$ and survival copula $C$.*

In particular, if $\mathbf{X}$ is a random vector with survival function $\bar{H}$ and continuous marginals $F_1, \ldots, F_d$ and $\mathbf{U}$ is a random vector distributed as the survival copula $C$ of $\bar{H}$, we have that

$$\mathbf{U} \stackrel{\mathrm{d}}{=} (\bar{F}_1(X_1), \ldots, \bar{F}_d(X_d)) \quad \text{and} \quad \mathbf{X} \stackrel{\mathrm{d}}{=} (\bar{F}_1^{-1}(U_1), \ldots, \bar{F}_d^{-1}(U_d)).$$

We are now in a position to turn our attention to Archimedean copulas. The latter were originally characterized by associativity and the property that $C(u, u) < u$ for any $u \in [0, 1]$; the present paper uses a more common definition based on the generator $\psi$.

DEFINITION 2.2. A nonincreasing and continuous function $\psi : [0, \infty) \to [0, 1]$ which satisfies the conditions $\psi(0) = 1$ and $\lim_{x \to \infty} \psi(x) = 0$ and is strictly decreasing on $[0, \inf\{x : \psi(x) = 0\})$ is called an *Archimedean generator*. A $d$-dimensional copula $C$ is called *Archimedean* if it permits the representation

$$C(\mathbf{u}) = \psi(\psi^{-1}(u_1) + \cdots + \psi^{-1}(u_d)), \qquad \mathbf{u} \in [0,1]^d$$

for some Archimedean generator $\psi$ and its inverse $\psi^{-1} : (0, 1] \to [0, \infty)$ where, by convention, $\psi(\infty) = 0$ and $\psi^{-1}(0) = \inf\{u : \psi(u) = 0\}$.

Note that several authors define Archimedean copulas in terms of $\psi^{-1}$ rather than $\psi$. The reason the above definition is favored here is that it leads, in this context, to simpler expressions, as will soon be apparent from the discussions below.

A closer look at (1.1) readily reveals that $\psi(\psi^{-1}(u_1) + \cdots + \psi^{-1}(u_d))$ always satisfies the boundary conditions (i) and (ii) of Definition 2.1. Therefore, an Archimedean generator $\psi$ defines a $d$-dimensional copula via (1.1) if and only if $\psi(\psi^{-1}(u_1) + \cdots + \psi^{-1}(u_d))$ is quasi-monotone. The following well-known result, which is Theorem 6.3.2 of [37], states conditions under which quasi-monotonicity holds in the bivariate case.

PROPOSITION 2.1. *Let $\psi$ be an Archimedean generator in the sense of Definition 2.2. Then the function given by*

$$\psi(\psi^{-1}(u) + \psi^{-1}(v))$$

*for $u, v \in [0,1]$ is a copula if and only if $\psi$ is convex.*

Proposition 2.1 does not extend to dimensions $d \geq 3$, as illustrated below.

EXAMPLE 2.1. It is a simple matter to check that the function $\psi(x) = \max(1-x, 0)$ is a convex Archimedean generator. However,

$$\psi(\psi^{-1}(u_1) + \cdots + \psi^{-1}(u_d)) = \max(u_1 + \cdots + u_d - d + 1, 0),$$

which is the Fréchet–Hoeffding lower bound $W(u_1, \ldots, u_d)$. As is well known, this function is not a copula for $d \geq 3$, which can be seen by noting that it assigns negative mass to $[1/2, 1]^d$.

Consequently, stronger requirements on $\psi$ are needed. As will be shown in the sequel, these are based on the notion of $d$-monotone functions introduced below.

DEFINITION 2.3. A real function $f$ is called *d-monotone* in $(a,b)$, where $a, b \in \overline{\mathbb{R}}$ and $d \geq 2$, if it is differentiable there up to the order $d-2$ and the derivatives satisfy

$$(-1)^k f^{(k)}(x) \geq 0, \qquad k = 0, 1, \ldots, d-2$$

for any $x \in (a,b)$ and further if $(-1)^{d-2} f^{(d-2)}$ is nonincreasing and convex in $(a,b)$. For $d = 1$, $f$ is called *1-monotone* in $(a,b)$ if it is nonnegative and nonincreasing there. If $f$ has derivatives of all orders in $(a,b)$ and if $(-1)^k f^{(k)}(x) \geq 0$ for any $x$ in $(a,b)$, then $f$ is called *completely monotone*.

REMARK 2.1. Note that $d$-monotonicity as defined here is different from the notion of $k$-monotonicity as used, for example, in approximation. More specifically, a function $f:(a,b) \to \mathbb{R}$ is called $k$-monotone (or $k$-convex in [33] and elsewhere) if its divided difference of order $k$ is nonnegative. Without going into more detail at this point, we highlight the difference between the two concepts for $k = 2$: $k$-monotone then means simply convex, while a 2-monotone function in our sense is nonnegative, nonincreasing and convex. In particular, therefore, a $k$-monotone function is not necessarily $k'$-monotone for $k' < k$ which contrasts with the above defined notion of $d$-monotonicity. Furthermore, it can be shown that, provided it exists, the $k$th-order derivative of a $k$-monotone function is nonnegative, whereas in our case the derivatives alternate in sign.

Definition 2.3 can be extended to functions on not necessarily open intervals.

DEFINITION 2.4. A real function $f$ on an interval $I \subseteq \overline{\mathbb{R}}$ is $d$-monotone (completely monotone) on $I$, $d \in \mathbb{N}$, if it is continuous there and if $f$ restricted to the interior $I^\circ$ of $I$ is $d$-monotone (completely monotone) on $I^\circ$.

The central importance of $d$-monotonicity in the present paper is based on the following result, which relates $d$-monotonicity to the existence of survival functions.

PROPOSITION 2.2. *Let $f$ be a real function on $[0, \infty)$, $p \in [0,1]$ and $\bar{H}$ be specified by*

$$\bar{H}(\mathbf{x}) = f(\|\max(\mathbf{x}, \mathbf{0})\|_1) + (1-p)\mathbf{1}\{\mathbf{x} < \mathbf{0}\}, \qquad \mathbf{x} \in \mathbb{R}^d.$$

*Then $\bar{H}$ is a survival function on $\mathbb{R}^d$ if and only if $f$ is a $d$-monotone function on $[0, \infty)$ satisfying the boundary conditions $\lim_{x \to \infty} f(x) = 0$ and $f(0) = p$.*

The proof, together with that of Theorem 2.2 below, is deferred to Appendix A.

THEOREM 2.2. *Let $\psi$ be an Archimedean generator. Then $C:[0,1]^d \to [0,1]$ given by*

$$C(\mathbf{u}) = \psi(\psi^{-1}(u_1) + \cdots + \psi^{-1}(u_d)), \qquad \mathbf{u} \in [0,1]^d$$

*is a $d$-dimensional copula if and only if $\psi$ is $d$-monotone on $[0, \infty)$.*

Note that Theorem 2.2 entails Proposition 2.1 as a special case. A further straightforward consequence is as follows:

COROLLARY 2.1. *Suppose $\psi$ is an Archimedean generator which has derivatives up to order $d$ on $(0, \infty)$. Then $\psi$ generates an Archimedean copula if and only if $(-1)^k \psi^{(k)}(x) \geq 0$ for $k = 1, \ldots, d$.*

Obviously, $d$-monotonicity of an Archimedean generator $\psi$ does not imply that $\psi$ is $k$-monotone for $k > d$. In other words, $d$-monotone Archimedean generators do not necessarily generate Archimedean copulas in dimensions higher than $d$. This observation motivates the following notation:

DEFINITION 2.5. $\Psi_d$ denotes the class of $d$-monotone Archimedean generators for $d \geq 2$. In addition, $\Psi_\infty$ stands for the class of Archimedean generators which can generate an Archimedean copula in any dimension $d \geq 2$.

On the other hand, if a function is $d$-monotone for some $d \geq 2$, it is also $k$-monotone for any $1 \leq k \leq d$. This easy consequence of Definition 2.2 means that

$$\Psi_2 \supseteq \Psi_3 \supseteq \cdots.$$

Before giving examples, which highlight in particular that $\Psi_d \setminus \Psi_{d+1}$ is nonempty, we note that the verification of the $d$-monotonicity of a generator $\psi$ can be considerably simplified by means of the following result, which is excerpted from [41], Theorem 4:

PROPOSITION 2.3. *Let $d \geq 2$ be an integer and $\psi$ be an Archimedean generator in the sense of Definition 2.2. Then $\psi$ is $d$-monotone on $[0, \infty)$ if and only if $(-1)^{d-2}\psi^{(d-2)}$ exists on $(0, \infty)$ and is nonnegative, nonincreasing and convex there.*

For a complete account of the theory we can embed the well-known characterization of $\Psi_\infty$ due to [20] in our findings.

PROPOSITION 2.4. *An Archimedean generator $\psi$ belongs to $\Psi_\infty$ if and only if it is completely monotone on $[0, \infty)$.*

EXAMPLE 2.2. Consider the generator $\psi_d^{\mathbf{L}}(x) = (1-x)_+^{d-1}$ for some value of $d \geq 2$. The $\mathbf{L}$ superscript anticipates a result in Section 4 in which this generator is shown to generate a lower bound for $d$-dimensional Archimedean copulas. We can verify that $\psi_d^{\mathbf{L}}$ is $d$-monotone by computing $(\psi_d^{\mathbf{L}})^{(d-2)}(x) = (-1)^{d-2}(d-1)!(1-x)_+$ and observing that $(-1)^{d-2}(\psi_d^{\mathbf{L}})^{(d-2)}$ is nonnegative, nonincreasing and convex on $(0, \infty)$. It is not, however, $(d+1)$-monotone, since $(\psi_d^{\mathbf{L}})^{(d-1)}$ does not exist at $x = 1$. Hence, $\psi_d^{\mathbf{L}} \in \Psi_d \setminus \Psi_{d+1}$.

EXAMPLE 2.3. The Clayton copula family has generator $\psi_\theta(x) = (1 + \theta x)_+^{-1/\theta}$. It is well known that for $\theta > 0$ this generator is completely monotone and can be used to construct a copula in any dimension. The case $\theta = 0$ should be understood as the limit (from either side) $\psi_0(x) = \lim_{\theta \to 0}(1 + \theta x)_+^{-1/\theta} = \exp(-x)$, which generates the independence copula in any dimension.

The interesting case is $\theta < 0$. When $\theta = -1/(d-1)$ for some integer $d \geq 2$, then $\psi_\theta(x) = \psi_d^{\mathbf{L}}(\theta x)$ where $\psi_d^{\mathbf{L}}$ is the generator in Example 2.2. The argument of that example holds and the multiplicative $\theta$ is unimportant;

MULTIVARIATE ARCHIMEDEAN COPULAS 9$\psi_\theta$ generates copulas up to dimension $d$ and, in fact, generates exactly the same copulas as $\psi_d^{\mathbf{L}}$. For the general negative case let $\alpha = -1/\theta > 0$, write $\tilde{\psi}_\alpha(x) = \psi_{-1/\alpha}(x)$ and observe first that $\tilde{\psi}_\alpha$ is convex or 2-monotone when $\alpha \geq 1$. Moreover, when $d \geq 3$ and $\alpha \geq d-1$, then

$$(-1)^{d-2}\tilde{\psi}_\alpha^{(d-2)}(x) = \alpha^{-(d-2)} \prod_{k=0}^{d-3}(\alpha - k)\left(1 - \frac{x}{\alpha}\right)_+^{\alpha-d+2}$$

exists and is nonnegative, nonincreasing and convex on $(0, \infty)$, whereas for other $\alpha$ values the derivative will either fail to exist everywhere or fail to have these properties; thus, by Proposition 2.3, the generator $\tilde{\psi}_\alpha$ is $d$-monotone if and only if $\alpha \geq d-1$. Summarizing, $\psi_\theta$ is $d$-monotone for a particular $d \geq 2$ if and only if $\theta \geq -1/(d-1)$.

**3. Archimedean copulas and $\ell_1$-norm symmetric distributions.** For the remainder of the paper, assume that $d \geq 2$. The beauty of Kimberling's characterization of $\Psi_\infty$ lies in its combination with the well-known Bernstein–Widder theorem (see, e.g., [40]). The latter states that an Archimedean generator is completely monotone on $[0, \infty)$ precisely when it is the Laplace transform of a nonnegative random variable, say $W$. In this case, as is well known, the $d$-dimensional Archimedean copula generated by $\psi$ is a survival copula of the survival function

$$\bar{H}(\mathbf{x}) = \psi(\|\max(\mathbf{x}, \mathbf{0})\|_1) = E(e^{-\|\max(\mathbf{x}, \mathbf{0})\|_1 W}) = E(e^{-W\sum_{i=1}^d \max(x_i, 0)}),$$

which is the survival function of the random vector $\mathbf{X} = \mathbf{Y}/W$, where $\mathbf{Y} = (Y_1, \ldots, Y_d)$ is a vector of i.i.d. exponential variables, independent of $W$. $\mathbf{X}$ follows a mixed exponential distribution or frailty model of the kind that is often used to model dependent lifetimes; see [23, 29].

However, it is possible to give another representation for $\mathbf{X}$. We can use the well-known facts that the random vector $\mathbf{S}_d = \mathbf{Y}/\|\mathbf{Y}\|_1$ has a uniform distribution on the $d$-dimensional simplex and $\mathbf{S}_d$ and $\|\mathbf{Y}\|_1$ are independent to write $\mathbf{X} = R\mathbf{S}_d$ where $R = \|\mathbf{Y}\|_1/W$ and $R$ is independent of $\mathbf{S}_d$. This shows that $\mathbf{X}$ follows a so-called $\ell_1$-norm symmetric distribution, a mixture of uniform distributions on simplices also known as a simplicially contoured distribution.

Archimedean generators that are not completely monotone cannot appear as survival copulas of random vectors following frailty models, but they can all appear as survival copulas of random vectors following $\ell_1$-norm symmetric distributions. These distributions, introduced in [7], are formally defined as follows:

DEFINITION 3.1. A random vector $\mathbf{X}$ on $\mathbb{R}_+^d = [0, \infty)^d$ follows an $\ell_1$-norm symmetric distribution if and only if there exists a nonnegative random variable $R$ independent of $\mathbf{S}_d$ where $\mathbf{S}_d$ is a random vector distributed uniformly on the unit simplex $\mathcal{S}_d$,
$$\mathcal{S}_d = \{\mathbf{x} \in \mathbb{R}^d_+ : \|\mathbf{x}\|_1 = 1\},$$
so that $\mathbf{X}$ permits the stochastic representation
$$\mathbf{X} \stackrel{\mathrm{d}}{=} R\mathbf{S}_d.$$
The random variable $R$ is referred to as the *radial part* of $\mathbf{X}$ and its distribution as the *radial distribution*.

The general relationship between Archimedean copulas and $\ell_1$-norm symmetric distributions is quite similar to the relationship between Archimedean copulas with completely monotone generators and frailty models outlined above. The difference is that the role of the Laplace transform is taken by another integral transform, which we call the Williamson $d$-transform and define as follows:

DEFINITION 3.2. *Let $X$ be a nonnegative random variable with distribution function $F$ and $d \geq 2$ an integer. The Williamson $d$-transform of $X$ is a real function on $[0, \infty)$ given by*
$$\mathfrak{W}_d F(x) = \int_{(x,\infty)} \left(1 - \frac{x}{t}\right)^{d-1} dF(t) = \begin{cases} E\left(1 - \dfrac{x}{X}\right)^{d-1}_+, & \text{if } x > 0, \\ 1 - F(0), & \text{if } x = 0. \end{cases}$$
*The class of functions which are Williamson $d$-transforms of nonnegative random variables will be denoted by $\mathcal{W}_d$.*

REMARK 3.1. The Williamson $d$-transform belongs to the general class of Mellin–Stieltjes convolutions and is closely related to the Cesàro means. For more details, refer to [2], pages 194 and 246.

It is clear that a Williamson $d$-transform always exists and is right-continuous at 0 by the dominated convergence theorem. Far less obvious is the following result by Williamson [41]:

PROPOSITION 3.1. *Let $d \geq 2$ be an arbitrary integer. Then*

(i) $\mathcal{W}_d$ *consists precisely of the real functions $f$ on $[0, \infty)$ which are $d$-monotone on $[0, \infty)$ and satisfy the boundary conditions $\lim_{x \to \infty} f(x) = 0$ as well as $f(0) = p$ for $p \in [0, 1]$.*

(ii) *The distribution of a nonnegative random variable is uniquely given by its Williamson $d$-transform. If $f = \mathfrak{W}_d F$, then for $x \in [0, \infty)$, $F(x) = \mathfrak{W}_d^{-1} f(x)$, where*
$$\mathfrak{W}_d^{-1} f(x) = 1 - \sum_{k=0}^{d-2} \frac{(-1)^k x^k f^{(k)}(x)}{k!} - \frac{(-1)^{d-1} x^{d-1} f_+^{(d-1)}(x)}{(d-1)!}.$$

PROOF. The proof requires a reformulation of Williamson's original result in terms of distribution functions. This is carried out in Appendix B. □

It may be noted that the proof of Proposition 3.1 indicates that if $f \in \mathcal{W}_d$, then $f(0) = p$ if and only if the corresponding distribution function $F$ satisfies $F(0) = 1 - p$ (i.e., has an atom at 0 of size $1-p$ when $p < 1$). Furthermore, since any $f \in \mathcal{W}_d$ is strictly decreasing on $[0, \inf\{x : f(x) = 0\})$, we can conclude the following:

COROLLARY 3.1. $\Psi_d$ consists of the Williamson d-transforms of distribution functions $F$ of nonnegative random variables satisfying $F(0) = 0$.

The role of the Williamson $d$-transform in the study of $\ell_1$-norm symmetric distributions is indicated in Proposition 3.2 below, which also recalls two further properties of $\ell_1$-norm symmetric distributions that will be used in this paper. These and many other results on $\ell_1$-norm symmetric distributions have been derived by K. T. Fang, B. Q. Fang and collaborators; see the monograph [8].

PROPOSITION 3.2. Let $\mathbf{X} \stackrel{d}{=} R\mathbf{S}_d$ be a random vector which follows an $\ell_1$-norm symmetric distribution with radial distribution function $F_R$. Then

(i) The survival function $\bar{H}$ of $\mathbf{X}$ is given by

$$(3.1) \quad \bar{H}(\mathbf{x}) = \mathfrak{W}_d F_R(\|\max(\mathbf{x}, \mathbf{0})\|_1) + F_R(0)\mathbf{1}\{\mathbf{x} < \mathbf{0}\}, \qquad \mathbf{x} \in \mathbb{R}^d.$$

(ii) The density $\mathbf{X}$ exists if and only if $R$ has a density. In that case, it is given by $h(\|\mathbf{x}\|_1) = \Gamma(d)\|\mathbf{x}\|^{1-d} f_R(\|\mathbf{x}\|_1)$ where $f_R$ denotes the density of $R$.

(iii) If $P(\mathbf{X} = \mathbf{0}) = 0$, then $R \stackrel{d}{=} \|\mathbf{X}\|_1$ and $\mathbf{S}_d \stackrel{d}{=} \mathbf{X}/\|\mathbf{X}\|_1$.

PROOF. See [8], Theorems 5.1, 5.4 and 5.5. □

The important implication of this proposition is that the survival functions of $\ell_1$-norm symmetric distributions reduce to $\mathfrak{W}_d F_R(\|\mathbf{x}\|_1)$, a simple function of $\|\mathbf{x}\|_1$, whenever $\mathbf{x} \in \mathbb{R}_+^d$. In Proposition 3.3 below we extend this result to show the converse, namely that whenever a survival function is a function of $\|\mathbf{x}\|_1$ on $\mathbf{x} \in \mathbb{R}_+^d$ it must be the survival function of an $\ell_1$-norm symmetric distribution, provided it places no mass on the boundary of $\mathbb{R}_+^d$, except possibly at the origin.

PROPOSITION 3.3. Let $\mathbf{X}$ be a random vector on $\mathbb{R}_+^d$. Then the following are equivalent:

  (i) **X** *has an $\ell_1$-norm symmetric distribution.*
  (ii) *There exists a real function $f$ on $[0, \infty)$ so that the joint survival function $\bar{H}$ of **X** satisfies, for any $\mathbf{x} \in \mathbb{R}^d$,*

$$\bar{H}(\mathbf{x}) = f(\|\max(\mathbf{x}, \mathbf{0})\|_1) + (1 - f(0))\mathbf{1}\{\mathbf{x} < \mathbf{0}\}. \tag{3.2}$$

PROOF. The fact that an $\ell_1$-norm symmetric distribution satisfies (ii) is an immediate consequence of Proposition 3.2. The reversed implication can be established as follows: First, note that if (ii) holds, $f$ must be $d$-monotone on $[0, \infty)$ and satisfy $f(0) \in [0, 1]$ as well as $\lim_{x \to \infty} f(x) = 0$ by Proposition 2.2. Proposition 3.1 implies that $f \in \mathcal{W}_d$. Let $F_R$ be the distribution function $F_R$ such that $f = \mathfrak{W}_d F_R$. Then $\bar{H}$ is precisely of the form (3.1) and **X** must have an $\ell_1$-norm symmetric distribution with radial distribution $F_R$. □

REMARK 3.2. Proposition 3.3 is needed to fully establish the connection between $\ell_1$-norm symmetric distributions and Archimedean copulas. In [7] a related result is proved that does not go quite far enough: the authors introduce a class $T_d$ of distributions on $\mathbb{R}_+^d$ whose survival functions are of the form $f(\|\mathbf{x}\|_1)$ for some function $f$ on $[0, \infty)$ and argue that a member of $T_d$ is $\ell_1$-norm symmetric if $f$ is $d$-monotone and $f(0) = 1$. What Proposition 3.3 shows is that $T_d$, in fact, coincides with the class of $\ell_1$-norm symmetric distributions as long as we simply exclude certain behavior on the boundary of $\mathbb{R}_+^d$.

We are now in a position to bring all the pieces of the argument together in the following, which is the main result of this section:

THEOREM 3.1. (i) *Let **X** have a $d$-dimensional $\ell_1$-norm symmetric distribution with radial distribution $F_R$ satisfying $F_R(0) = 0$. Then **X** has an Archimedean survival copula with generator $\psi = \mathfrak{W}_d F_R$.*

(ii) *Let **U** be distributed according to the $d$-dimensional Archimedean copula $C$ with generator $\psi$. Then $(\psi^{-1}(U_1), \ldots, \psi^{-1}(U_d))$ has an $\ell_1$-norm symmetric distribution with survival copula $C$ and radial distribution $F_R$ satisfying $F_R = \mathfrak{W}_d^{-1} \psi$, that is,*

$$F_R(x) = 1 - \sum_{k=0}^{d-2} \frac{(-1)^k x^k \psi^{(k)}(x)}{k!} - \frac{(-1)^{d-1} x^{d-1} \psi_+^{(d-1)}(x)}{(d-1)!}, \tag{3.3}$$

$$x \in [0, \infty).$$

PROOF. (i) From Proposition 3.2, part (i), it follows that the survival function of **X** satisfies $\bar{H}(\mathbf{x}) = \psi(\|\max(\mathbf{x}, \mathbf{0})\|_1)$ on $\mathbb{R}^d$ where $\psi = \mathfrak{W}_d F_R$.

By Corollary 3.1 we know that $\psi \in \Psi_d$. As shown in the proof of Theorem 2.2, since the marginal survival functions of $\bar{H}$ are the continuous functions $\bar{F}_i(x) = \psi(\max(x,0))$, we can conclude that the unique survival copula of $\mathbf{X}$ is the Archimedean copula generated by $\psi$.

(ii) The survival function of $(\psi^{-1}(U_1), \ldots, \psi^{-1}(U_d))$ is given by $\bar{H}(x_1, \ldots, x_d) = C(\psi(x_1), \ldots, \psi(x_d)) = \psi(\|\max(\mathbf{x}, \mathbf{0})\|_1)$ on $\mathbb{R}^d$. From Proposition 3.3 we conclude that this is the survival function of an $\ell_1$-norm symmetric distribution and $C$ is its copula. Since $\psi \in \Psi_d$, it must be the Williamson $d$-transform of some distribution function $F_R$ of a nonnegative random variable $R$ satisfying $F_R(0) = 0$ and this $R$ is the radial part in the representation $(\psi^{-1}(U_1), \ldots, \psi^{-1}(U_d)) \stackrel{\mathrm{d}}{=} R\mathbf{S}_d$. Finally, $F_R = \mathfrak{W}_d^{-1}\psi$ follows by (ii) of Proposition 3.1. □

REMARK 3.3. If we consider radial random variables with point mass at the origin and use these to construct $\ell_1$-norm symmetric distributions, then Proposition 3.2, part (i), shows that these continue to have survival functions determined by the Williamson $d$-transform $f(x) = \mathfrak{W}_d F_R(x)$, albeit in a case where $f(0) < 1$ and $f \notin \Psi_d$. The $\ell_1$-norm symmetric distribution does not have a unique copula since it has discontinuous univariate margins. We enter the treacherous world of copulas for discrete random variables; see [13].

In light of the second part of Theorem 3.1 it will be useful to define the notion of an $\ell_1$-norm symmetric distribution associated with a particular Archimedean generator.

DEFINITION 3.3. Let $\psi \in \Psi_d$ and let $F_R = \mathfrak{W}_d^{-1}\psi$ be the inverse Williamson $d$-transform of $\psi$ specified by (3.3). Then $F_R$ will be known as the radial distribution associated with $\psi$ in dimension $d$. If $R \sim F_R$ is a random variable independent of $\mathbf{S}_d$, where $\mathbf{S}_d$ is uniformly distributed on $\mathcal{S}_d$, the distribution of $R\mathbf{S}_d$ will be known as the $\ell_1$-norm symmetric distribution associated with $\psi$ in dimension $d$.

REMARK 3.4. Although there is a one-to-one relationship between $\psi$ and a particular radial distribution $F_R$ and corresponding $\ell_1$-norm symmetric distribution in a particular dimension, there is not a one-to-one relationship between the $d$-dimensional Archimedean copulas and the $\ell_1$-norm symmetric distributions. Consider the radial variables $R$ and $\tilde{R} = kR$ for $k > 0$. These have radial distribution functions $F_R(r)$ and $F_{\tilde{R}}(r) = F_R(r/k)$, respectively, and give rise to different $\ell_1$-norm symmetric distributions. Their Williamson $d$-transforms are $\psi(x) = \mathfrak{W}_d F_R(x)$ and $\tilde{\psi}(x) = \psi(x/k)$. However, these generate exactly the same Archimedean copula in dimension $d$.

We close this section by giving examples of the use of Theorem 3.1 which show how we can calculate the radial distribution associated with a given Archimedean copula generator.

EXAMPLE 3.1. Consider the generator $\psi_d^{\mathbf{L}}(x) = (1-x)_+^{d-1}$ of the copula $C_d^{\mathbf{L}}$ introduced in Example 2.2. Formula (3.3) easily yields that the radial distribution associated with $\psi_d^{\mathbf{L}}$ is, in fact, degenerate with $R = 1$ almost surely. Let $F_{\tilde{R}}$ be the distribution function of another random variable $\tilde{R}$ satisfying $\tilde{R} = k$ almost surely, where $k > 0$. Theorem 3.1, part (i), implies that the $d$-dimensional $\ell_1$-norm symmetric distribution with radial distribution $F_{\tilde{R}}$ has a survival copula with generator given by $\psi(x) = \mathfrak{W}_d F_{\tilde{R}}(x) = (1 - x/k)_+^{d-1}$. However, for any $k > 0$, this generates the same copula as $\psi_d^{\mathbf{L}}(x) = (1-x)_+^{d-1}$.

EXAMPLE 3.2. As is well known, the independence copula $\Pi(u_1, \ldots, u_d)$ is Archimedean with generator $\psi(x) = \exp(-x)$. It is again easy to show by means of (3.3) that the radial distribution associated with this generator is the Erlang distribution (or, alternatively, gamma distribution with parameter $d$) given by its density

$$f_R(x) = e^{-x} \frac{x^{d-1}}{(d-1)!}, \qquad x \in [0, \infty).$$

The survival function of the associated $\ell_1$-norm symmetric distribution is $\bar{H}(\mathbf{x}) = e^{-\sum_{i=1}^d x_i}$ by (3.1), which is obviously the survival function of independent standard exponential variables. Note that, in conjunction with Proposition 3.2, part (iii), this shows that a simple way of generating a random vector $\mathbf{S}_d$ on the unit simplex is to generate $d$ independent standard exponentials and divide each by the sum of the $d$ variables.

EXAMPLE 3.3. Consider the Clayton copula family of Example 2.3 with generators $\psi_\theta(x) = (1 + \theta x)_+^{-1/\theta}$. Assume that a particular $d \geq 2$ and a particular $\theta \geq -1/(d-1)$ are given. Suppose first that $\theta < 0$, set $\alpha = -1/\theta > 0$ and consider the generator $\tilde{\psi}_\alpha(x) = \psi_{-1/\alpha}(x)$. Using (3.3) the associated radial distribution may be calculated to be

$$\begin{aligned}
F_R(x) &= 1 - \sum_{k=0}^{d-1} \frac{(-1)^k x^k \tilde{\psi}_\alpha^{(k)}(x)}{k!} \\
&= 1 - \sum_{k=0}^{d-1} \frac{x^k}{k! \alpha^k} \prod_{j=0}^{k-1} (\alpha - j) \left(1 - \frac{x}{\alpha}\right)^{\alpha-k} \\
&= 1 - \sum_{k=0}^{d-1} \binom{\alpha}{k} \left(\frac{x}{\alpha}\right)^k \left(1 - \frac{x}{\alpha}\right)^{\alpha-k},
\end{aligned}$$

(3.4)

where $\binom{\alpha}{k}$ denotes the extended binomial coefficient. Moreover, $F_R(x) = 0$ and $F_R(x) = 1$ if, respectively, $x < 0$ and $x \geq \alpha$. Observe first that, if $\alpha = d-1$, then this distribution simplifies to point mass at $x = d-1$ as in Example 3.1, which is as we would expect. It is also an easy limiting calculation to establish that, as $\alpha \to \infty$, the distribution function of $R$ converges pointwise to the Erlang distribution in Example 3.2.

If $\theta > 0$, we can proceed in a similar manner by differentiating $\psi_\theta(x)$ to show that for $x \in [0, \infty)$ we have

$$(3.5) \qquad F_R(x) = 1 - \sum_{k=0}^{d-1} \frac{\prod_{j=0}^{k-1}(1+j\theta)}{k!} x^k (1+\theta x)^{-(1/\theta+k)}.$$

Again, as $\theta \to 0$ this converges pointwise to the Erlang distribution function in Example 3.2.

For illustration, consider the case when $d = 3$. If $\theta = -0.5$, the distribution of $R$ is point mass at $x = 2$. For $\theta > -0.5$, $R$ has a density and the copula is absolutely continuous. The radial densities are given in Figure 1 for the cases when $\theta$ takes the values $-0.3334$, $-0.3$, $-0.2$, $0$ and $0.2$.

**4. A new perspective on Archimedean copulas.** The close relationship between Archimedean copulas and $\ell_1$-norm symmetric distributions delineated in Theorem 3.1 can be exploited in various ways. The important ingredient is the fact that the distribution $F_R$ of the radial part of an $\ell_1$-norm symmetric distribution can be *explicitly* retrieved from $\mathfrak{W}_d F_R$ and hence from $\bar{H}$ by means of (3.3), the inversion formula for Williamson $d$-transforms. From this point of view, working with Williamson $d$-transforms is more convenient than working with Laplace transforms for which the inversion formula cannot always be evaluated explicitly.

The simple stochastic structure of $\ell_1$-norm symmetric distributions can help us to explore analytical as well as dependence properties of Archimedean copulas and several results of this kind are given in this section.

4.1. *Singular components of d-dimensional Archimedean copulas.* The well-known Lebesgue decomposition theorem yields that any copula $C$ can be written as

$$C = C_A + C_S,$$

where $C_A$ and $C_S$ are, respectively, the absolutely continuous and singular components of $C$. In other words, $C_A$ is a distribution function of a (finite) measure on $\mathbb{R}^d$ which is absolutely continuous w.r.t. the Lebesgue measure on $\mathbb{R}^d$, that is, it has Lebesgue density. On the other hand, the measure

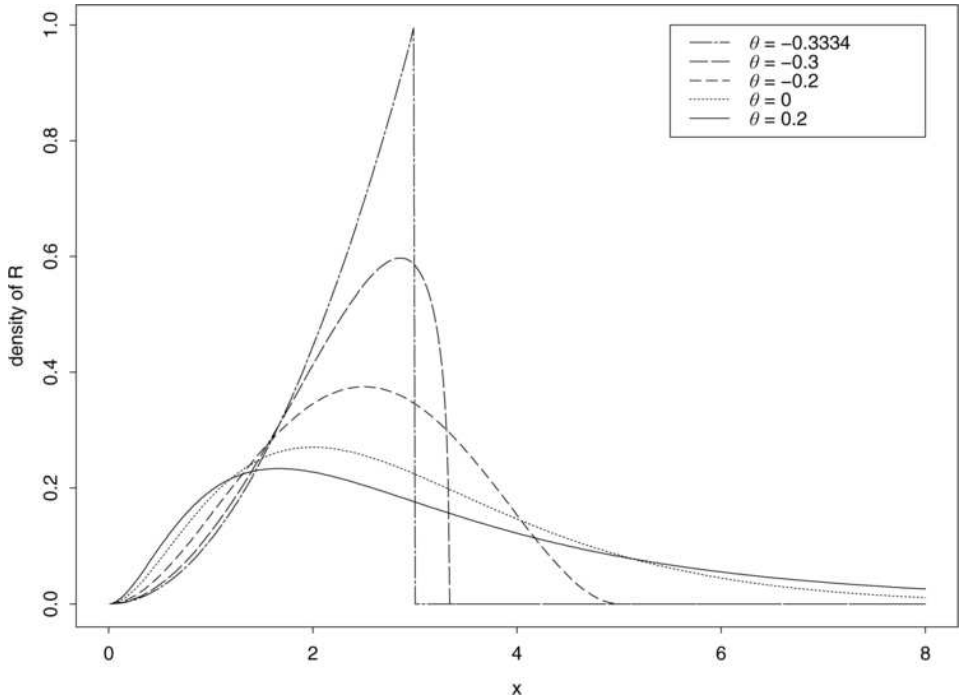

FIG. 1. *Radial densities for the trivariate Clayton copula when $\theta$ equals $-0.3334$, $-0.3$, $-0.2$, $0$ and $0.2$.*

induced by $C_S$ is singular, that is, it is concentrated on a set of Lebesgue measure zero.

Generally, the study of either component of a copula is not easy. However, in view of maximum likelihood estimation, it may be convenient to know whether a copula $C$ has a density. In this section, we obtain more concrete results when $C$ is Archimedean. Furthermore, singular components of Archimedean copulas in the bivariate case are discussed in [12] and [27]; we extend these results to higher dimensions.

To begin with, we should point out one potential fallacy. Let **X** denote a random vector with survival function $\bar{H}$ and continuous marginal survival functions $\bar{F}_1, \ldots, \bar{F}_d$. Then, according to Sklar's theorem, there exists a unique survival copula $C$ of $\bar{H}$. Furthermore, as mentioned in Section 2, if **U** is a random vector whose distribution function is $C$, then

$$\mathbf{U} \stackrel{\mathrm{d}}{=} (\bar{F}_1(X_1), \ldots, \bar{F}_d(X_d)) \quad \text{and} \quad \mathbf{X} \stackrel{\mathrm{d}}{=} (\bar{F}_1^{-1}(U_1), \ldots, \bar{F}_d^{-1}(U_d)).$$

In other words, if $P^C$ and $P^H$ denote, respectively, the probability measures induced by $C$ and $\bar{H}$, then $P^C$ is an image measure of $P^H$ with respect to a certain transformation and vice versa. Therefore, it may be tempting to believe that if $P^H$ is absolutely continuous, then so is $P^C$ and vice versa. This may be true when $\bar{F}_i$ and $\bar{F}_i^{-1}$ obey the usual smoothness conditions required by the classical change of variable formula for Lebesgue integrals. In general, however, it is only guaranteed that $\bar{F}_i$ is differentiable almost everywhere and to investigate the relationship between the continuity properties of $P^C$ and $P^H$ requires more involved change of variable formulas, which goes beyond the scope of this paper. Readers who are more interested in this subject should consult [17] for further information and references. We do provide a simple counterexample, however.

EXAMPLE 4.1. Take $F_{\mathfrak{C}}$ to be the Cantor function (see [16], Example 18.8). As is well known, $F_{\mathfrak{C}}:[0,1] \to [0,1]$ is strictly increasing, continuous and satisfies $F_{\mathfrak{C}}(0) = 0$ as well as $F_{\mathfrak{C}}(1) = 1$ but $F'_{\mathfrak{C}}(t) = 0$ a.e. in $[0,1]$. In particular, therefore, $F_{\mathfrak{C}}$ provides an example of a singular univariate distribution function. Now, consider the following survival function $\bar{H}$ on $\mathbb{R}^d$:

$$\bar{H}(\mathbf{x}) = \prod_{i=1}^{d} \bar{F}_{\mathfrak{C}}(x_i), \qquad \mathbf{x} \in \mathbb{R}^d.$$

It is immediately clear that $\bar{H}$ does not induce an absolutely continuous probability measure. On the other hand, the marginals $\bar{F}_{\mathfrak{C}}$ of $\bar{H}$ are continuous and there exists a unique survival copula corresponding to $\bar{H}$ by Sklar's theorem. In fact, the latter is the independence copula $\Pi$ given by $\Pi(\mathbf{u}) = \prod_{i=1}^{d} u_i$. Clearly, $\Pi$ has a density.

In the case of Archimedean copulas, however, absolute continuity can be led back to absolute continuity of the associated $\ell_1$-norm symmetric distribution as follows:

PROPOSITION 4.1. *Let $C$ be a d-dimensional Archimedean copula with generator $\psi \in \Psi_d$. Let $H$ stand for the distribution function of the $\ell_1$-norm symmetric distribution associated with $\psi$. Then*

(i) *$C$ is absolutely continuous if and only if $H$ is.*
(ii) *If $\psi \in \Psi_{d+1}$, then $C$ is absolutely continuous.*

PROOF. Refer to Appendix C. □

It may be noted that Proposition 4.1 implies in particular that all lower-dimensional margins of an Archimedean copula are absolutely continuous.

Proposition 3.2 states that an $\ell_1$-norm symmetric distribution is absolutely continuous if and only if its radial part is. Furthermore, the distribution function of the radial part is retrievable from $\psi$ using the inversion (3.3). These observations yield the following:

PROPOSITION 4.2. *Let $C$ be a $d$-dimensional Archimedean copula with generator $\psi$. Then $C$ is absolutely continuous if and only if $\psi^{(d-1)}$ exists and is absolutely continuous in $(0, \infty)$. Furthermore, if $C$ is absolutely continuous, then its density $c$ is given by*

$$c(\mathbf{u}) = \frac{\psi^{(d)}(\psi^{-1}(u_1) + \cdots + \psi^{-1}(u_d))}{\psi'(\psi^{-1}(u_1)) \cdots \psi'(\psi^{-1}(u_1))}$$

*for almost all $\mathbf{u} \in (0,1)^d$.*

PROOF. Verification of this claim may again be found in Appendix C. □

Note that, in particular, the Archimedean generators considered in Corollary 2.1 obey the conditions of Proposition 4.2.

Let us now summarize where we stand. If $C$ is a $d$-dimensional Archimedean copula, its generator $\psi$ belongs to either $\Psi_{d+1}$ or $\Psi_d \setminus \Psi_{d+1}$. In the former case, $C$ has a density by means of Proposition 4.1. In the latter, however, the density of $C$ may or may not exist. Examples below illustrate situations that arise for $\psi \in \Psi_d \setminus \Psi_{d+1}$.

EXAMPLE 4.2. $\psi^{(d-1)}$ *is absolutely continuous on $(0, \infty)$.* In this case, Proposition 4.2 ensures that $C$ has a density. An example in the bivariate case is as follows: Take $\psi$ on $[0, \infty)$ given by

$$\psi(x) = \begin{cases} \frac{8}{7}x^2 - \frac{16}{7}x + 1, & \text{if } x \in [0, 0.5), \\ \frac{16}{7}x^2 - \frac{24}{7}x + \frac{9}{7}, & \text{if } x \in [0.5, 0.75), \\ 0, & \text{if } x \in [0.75, \infty). \end{cases}$$

Because $\psi'(x) = \frac{16}{7}(x-1)$ for $x \in [0, 0.5)$ and $\psi'(x) = \frac{16}{7}(2x - \frac{3}{2})$ for $x \in [0.5, 0.75)$, $-\psi'$ is nonincreasing and nonnegative on $(0, \infty)$. In particular, therefore, $\psi \in \Psi_2$. Furthermore, $\psi'$ is absolutely continuous in $[0, \infty)$ and the bivariate Archimedean copula induced by $\psi$ is absolutely continuous. On the other hand, however, it can be easily verified that $-\psi'$ is not convex. Consequently, $\psi$ is not 3-monotone on $(0, \infty)$ and hence not in $\Psi_3$.

EXAMPLE 4.3. $\psi^{(d-1)}$ *is continuous but not absolutely continuous on $(0, \infty)$.* We again provide an example in the bivariate case. Consider the function

$$\psi(x) = \begin{cases} 1 - \frac{1}{c} \int_0^x F_{\mathfrak{C}}(1-t)\, dt, & \text{if } t \in [0,1], \\ 0, & \text{if } t > 1, \end{cases}$$

where $c = \int_0^1 F_{\mathfrak{C}}(1-t)\,dt$ and $F_{\mathfrak{C}}$ is the Cantor function introduced in Example 4.1. It is first clear that $\psi(0) = 1$, $\psi(0) = 0$ and that $\psi$ is continuous on $[0, \infty)$. Furthermore,

$$\psi'(x) = \begin{cases} -\dfrac{1}{c}F_{\mathfrak{C}}(1-x), & \text{if } t \in [0,1], \\ 0, & \text{if } t > 1, \end{cases}$$

for any $x \in [0, \infty)$. In particular, $-\psi'$ is continuous, nonnegative and nonincreasing and thus $\psi \in \Psi_2$. However, $\psi'$ is singular on $[0,1]$. It may be noted that the radial part of the associated $\ell_1$-norm symmetric distribution is also purely singular because, for almost all $x \in [0,1]$, $F_R'(x) = \frac{1}{c}F_{\mathfrak{C}}'(1-x) = 0$.

EXAMPLE 4.4. $\psi_+^{(d-1)}$ *has jumps.* Situations of this kind arise easily upon considering radial parts with atoms. To provide an example, set $R$ to be a random variable independent of $\mathbf{S}_d$ which follows the geometric distribution, that is, $P(R = i) = p(1-p)^{i-1}$ for $i \in \mathbb{N}$. The random vector $\mathbf{X} \stackrel{\mathrm{d}}{=} R\mathbf{S}_d$ is then $\ell_1$-norm symmetrically distributed and its survival copula $C$ is Archimedean with generator

$$\psi(x) = \mathfrak{W}_d F_R(x) = \sum_{i=1}^{\infty}\left(1 - \frac{x}{i}\right)_+^{d-1} p(1-p)^{i-1}, \qquad x \in [0,\infty).$$

In particular, because $F_R$ has jumps and $\psi$ and $\psi^{(k)}$, $k = 1,\ldots,d-2$ are continuous, (3.3) implies that $\psi_+^{(d-1)}$ has jumps and is thus not absolutely continuous. This confirms that $C$ does not have density, which is already clear from the fact that $R$ has atoms.

So far, we have used the relationship between $\ell_1$-norm symmetric distributions and Archimedean copulas in order to examine the existence of densities. As the rest of this section shows, the latter can be successfully employed for the investigation of the singular part of Archimedean copulas as well. We begin with a simple example:

EXAMPLE 4.5. A similar situation to Example 4.4 arises when we consider $d = 2$ and a radial part satisfying $P(R = 1) = 2/3$ and $P(R = 2) = 1/3$. The Archimedean generator of the survival copula of $R\mathbf{S}_2$ is then

$$\psi(x) = \mathfrak{W}_d F_R(x) = \frac{2}{3}(1-x)_+ + \frac{1}{3}\left(1-\frac{x}{2}\right)_+.$$

By the same argument as in Example 4.4, the bivariate Archimedean copula $C$ induced by $\psi$ does not have a density. Figure 2 illustrates that the $\ell_1$-norm symmetric distribution with radial part $R$ is concentrated on two simplices,

$$\mathcal{S}_2(1) = \{\mathbf{x} \in \mathbb{R}_+^2 : \|\mathbf{x}\|_1 = 1\} \quad \text{and} \quad \mathcal{S}_2(2) = \{\mathbf{x} \in \mathbb{R}_+^2 : \|\mathbf{x}\|_1 = 2\}.$$

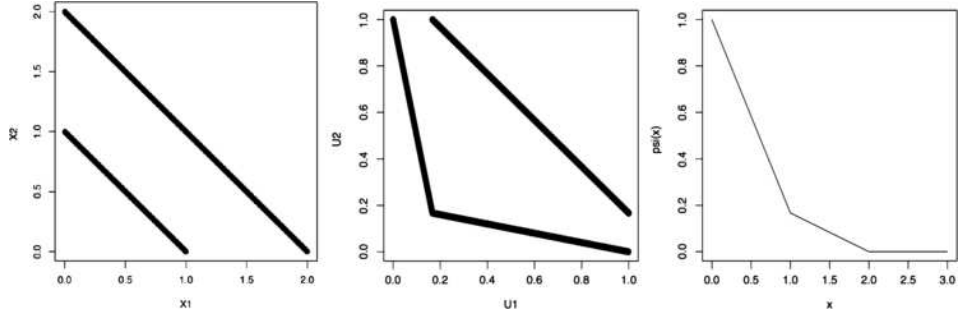

FIG. 2. *10,000 simulated values from the bivariate $\ell_1$-norm symmetric distribution of Example 4.5 (left panel). The middle panel shows the corresponding survival copula and the right panel shows its generator $\psi$.*

$C$ is also purely singular as can be seen in the middle panel of Figure 2. It is not difficult to verify that the support of $C$ is $\mathcal{A}_1 \cup \mathcal{A}_2$, where

$$\mathcal{A}_i = \{(u_1, u_2) \in [0,1]^2 : \psi^{-1}(u_1) + \psi^{-1}(u_2) = i\}, \qquad i = 1, 2.$$

Note that if $\mathbf{U} \sim C$, then $P(\mathbf{U} \in \mathcal{A}_1) = 2/3$ and $P(\mathbf{U} \in \mathcal{A}_2) = 1/3$.

The findings of Example 4.5 can be generalized. To do so, we first need to introduce additional notation. A *level set* $L(s)$ of a copula $C$ is given by

$$L(s) = \{\mathbf{u} \in [0,1]^d : C(\mathbf{u}) = s\}$$

for $s \in [0, 1]$. For a $d$-dimensional Archimedean copula the level sets take the form

$$L(s) = \begin{cases} \left\{\mathbf{u} \in [0,1]^d : \sum_{i=1}^d \psi^{-1}(u_i) = \psi^{-1}(s)\right\}, & \text{if } s \in (0, 1], \\ \left\{\mathbf{u} \in [0,1]^d : \sum_{i=1}^d \psi^{-1}(u_i) \geq \psi^{-1}(0)\right\}, & \text{if } s = 0. \end{cases}$$

Proposition 4.3 below, which is a high-dimensional extension of the results due to [12] and Alsina, Frank and Schweizer (see [27], Section 4.3), indicates the mass placed by an Archimedean copula on its level sets.

PROPOSITION 4.3. *Let $C$ be a $d$-dimensional Archimedean copula with generator $\psi$ and let $P^C$ stand for the probability measure on $[0,1]^d$ induced by $C$. Then*

$$P^C(L(s)) = \frac{(-1)^{d-1}(\psi^{-1}(s))^{d-1}}{(d-1)!}(\psi_-^{(d-1)}(\psi^{-1}(s)) - \psi_+^{(d-1)}(\psi^{-1}(s)))$$

*for $s \in (0,1]$. Furthermore, if $\psi^{-1}(0) = \infty$, then $P^C(L(0)) = 0$. Otherwise,*

$$P_C(L(0)) = P^C\left\{\mathbf{u} \in [0,1]^d : \sum_{i=1}^d \psi^{-1}(u_i) = \psi^{-1}(0)\right\}$$
$$= \frac{(-1)^{d-1}(\psi^{-1}(0))^{d-1}}{(d-1)!}\psi_-^{(d-1)}(\psi^{-1}(0)).$$

PROOF. Refer to Appendix C. □

4.2. *The Kendall function in $d$ dimensions.* For a bivariate Archimedean copula $C$ and a random vector $(U_1, U_2) \stackrel{\mathrm{d}}{=} C$, the Kendall function $K_C$ (also known as the bivariate probability integral transform) defined as the distribution function of $C(U_1, U_2)$ is a cornerstone of nonparametric inference for Archimedean copulas as introduced in [15]. Although this paper is not devoted to estimation, it is nonetheless worth noting that Theorem 3.1 allows an easy generalization of several of the findings about $K_C$ to higher dimensions. Results from this section also nicely complement those in [1], Example 3.

Let $C$ be a $d$-dimensional Archimedean copula with generator $\psi$ and $\mathbf{U} \sim C$ be a random vector. The function $K_C$ is, in analogy to the bivariate case, given as the distribution function of the random variable $C(\mathbf{U})$, that is,

$$K_C(x) = P(C(U_1, \ldots, U_d) \leq x).$$

Clearly, $C(\mathbf{U})$ is concentrated on $[0,1]$. For bivariate Archimedean copulas, $K_C(x)$ can be given explicitly in terms of the generator $\psi$ and its left-hand derivative (see [15], Proposition 1.1). To establish an analogous result in the $d$-dimensional case, observe first that

$$P(C(U_1, \ldots, U_d) \leq x) = P\left(\sum_{i=1}^d \psi^{-1}(U_i) \geq \psi^{-1}(x)\right).$$

Furthermore, let $\mathbf{X}$ be a random vector which follows the $\ell_1$-norm symmetric distribution associated with $\psi$. Sklar's theorem then implies

$$\mathbf{X} \stackrel{\mathrm{d}}{=} (\psi^{-1}(U_1), \ldots, \psi^{-1}(U_d)).$$

In particular, Proposition 3.2 yields that the radial part $R$ of $\mathbf{X}$ satisfies

$$R \stackrel{\mathrm{d}}{=} \|\mathbf{X}\|_1 \stackrel{\mathrm{d}}{=} \sum_{i=1}^d \psi^{-1}(U_i)$$

and is independent of $\mathbf{X}/\|\mathbf{X}\|_1$. The consequences are now immediate.

PROPOSITION 4.4. *Let $C$ be a $d$-dimensional Archimedean copula with generator $\psi$ and let $\mathbf{U}$ be a random vector with $\mathbf{U} \sim C$. Then*

$$\sum_{i=1}^{d} \psi^{-1}(U_i) \quad and \quad \left( \frac{\psi^{-1}(U_1)}{\sum_{i=1}^{d} \psi^{-1}(U_i)}, \ldots, \frac{\psi^{-1}(U_d)}{\sum_{i=1}^{d} \psi^{-1}(U_i)} \right)$$

*are independent. Moreover, for any $j = 1, \ldots, d$, the random variable*

$$V_j = \left( 1 - \frac{\psi^{-1}(U_j)}{\sum_{i=1}^{d} \psi^{-1}(U_i)} \right)^{d-1}$$

*is standard uniform.*

PROOF. What remains to be shown is the uniformity of $V_j$, $j = 1, \ldots, d$. This follows easily by observing the fact that the univariate marginals of the uniform distribution on the simplex $\mathcal{S}_d$ follow the $\mathrm{Beta}(1, d-1)$ distribution, that is, their survival functions are given by $\psi_d^{\mathbf{L}}(x) = (1-x)_+^{d-1}$. □

PROPOSITION 4.5. *Let $C$ be a $d$-dimensional Archimedean copula generated by $\psi$. Then*

$$K_C(x) = \begin{cases} \dfrac{(-1)^{d-1}(\psi^{-1}(0))^{d-1}}{(d-1)!} \psi_-^{(d-1)}(\psi^{-1}(0)), & \text{if } x = 0, \\ \displaystyle\sum_{k=0}^{d-2} \dfrac{(-1)^k (\psi^{-1}(x))^k \psi^{(k)}(\psi^{-1}(x))}{k!} \\ \quad + \dfrac{(-1)^{d-1}(\psi^{-1}(x))^{d-1} \psi_-^{(d-1)}(\psi^{-1}(x))}{(d-1)!}, & \text{if } x \in (0,1]. \end{cases}$$

PROOF. Because $K_C(x) = P(R \geq \psi^{-1}(x))$ where $R$ denotes the radial part of the $\ell_1$-norm symmetric distribution associated with $\psi$, the assertion follows directly from part (ii) of Theorem 3.1 and from Proposition 4.3. □

4.3. *Archimedean copulas are bounded below with respect to the PLOD ordering.* Recall that the positive lower orthant dependence (PLOD) ordering is a partial ordering on the set of all $d$-dimensional copulas given as follows [19]:

$$C_1 \prec_{\mathrm{cL}} C_2 \quad \Leftrightarrow \quad C_1(\mathbf{u}) \leq C_2(\mathbf{u}) \qquad \text{for any } \mathbf{u} \in [0,1]^d.$$

It is a well-established fact that any bivariate copula $C$ satisfies $W \prec_{\mathrm{cL}} C$ where $W$ is the Fréchet–Hoeffding lower bound specified in Example 2.1. In dimensions greater than two, it still holds that $W \prec_{\mathrm{cL}} C$, with the unfortunate exception that $W$ is no longer a copula. Furthermore, it can be shown that for $d \geq 3$, no sharp lower bound of the set of all copulas with respect to the PLOD ordering exists (see [27], Theorem 2.10.13).

The situation for Archimedean copulas is different, however. As is well known, an Archimedean copula $C$ whose generator is completely monotone is positive lower orthant dependent, that is, $\Pi \prec_{\mathrm{cL}} C$ where $\Pi$ is the independence copula given by $\Pi(\mathbf{u}) = \prod_{i=1}^{d} u_i$ (see [27], Corollary 4.6.3.). Proposition 4.6 below shows that a similar statement can also be made in the general case.

PROPOSITION 4.6. *Let $C$ be a $d$-dimensional Archimedean copula. Then*

$$C_d^{\mathbf{L}} \prec_{\mathrm{cL}} C,$$

*where $C_d^{\mathbf{L}}$ is the $d$-dimensional Archimedean copula with generator $\psi_d^{\mathbf{L}}$ of Example 2.2.*

PROOF. Refer to Appendix C. □

Figure 3 below illustrates the lower bound $C_d^{\mathbf{L}}$ in dimension $d=3$.

A useful implication of Proposition 4.6 is that if $C$ is the bivariate marginal distribution function of a $d$-dimensional Archimedean copula with generator $\psi \in \Psi_d$ and $\rho$ an arbitrary measure of concordance in the sense of [35], then

$$\rho(C_d^{\mathbf{L}}) \leq \rho(C),$$

where $\rho(C_d^{\mathbf{L}})$ refers to $\rho$ applied to the bivariate marginal distribution of $C_d^{\mathbf{L}}$. This applies in particular to Kendall's $\tau$, which is a concordance measure given by

$$\tau(C) = 4 \int_0^1 \int_0^1 C(u_1, u_2) \, dC(u_1, u_2) - 1 = 4E(C(U_1, U_2)) - 1,$$

where $(U_1, U_2)$ refers to a random vector distributed according to $C$.

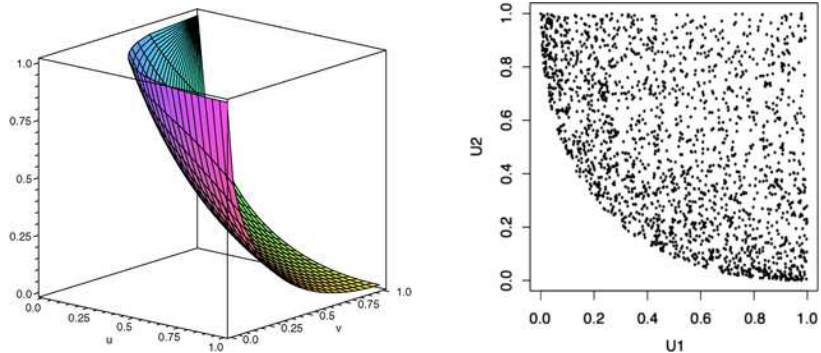

FIG. 3. *The left panel shows the support of the lower bound Archimedean copula $C_d^{\mathbf{L}}$ for $d=3$. Two thousand simulated points from the bivariate margin of $C_3^{\mathbf{L}}$ are depicted in the right panel.*

For bivariate Archimedean copulas, it was discovered in [11] that $\tau(C)$ can be expressed in terms of the corresponding Archimedean generator,

$$\tau(C) = 4 \int_0^1 \psi^{-1}(t)\psi'(\psi^{-1}(t))\,dt + 1$$

(4.1)

$$= 1 - 4 \int_0^{\psi^{-1}(0)} t[\psi'(t)]^2\,dt,$$

where the second expression results from the first by a change of variable. The following simple observation offers a geometric interpretation in terms of the associated $\ell_1$-norm symmetric distribution:

PROPOSITION 4.7. *Let $C$ be a bivariate Archimedean copula with generator $\psi$ and $R$ be the radial part of the $\ell_1$-norm symmetric distribution corresponding to $C$ in the sense of Definition 3.3. Then*

$$\tau(C) = 4E(\psi(R)) - 1.$$

PROOF. Let $(U_1, U_2)$ be a random vector distributed as $C$. The assertion is then an immediate consequence of the fact that $\psi(\psi^{-1}(U_1) + \psi^{-1}(U_2)) \stackrel{d}{=} \psi(R)$, which has been established in Section 4.2. □

The lower bound on Kendall's tau $\tau(C)$ is now straightforward.

COROLLARY 4.1. *Let $C$ be the bivariate margin of a $d$-dimensional Archimedean copula with generator $\psi$. Then $\tau(C) \geq -1/(2d-3)$.*

PROOF. The proof follows from (4.1) using elementary calculus. □

**5. Practical implications.** The results and connections established in this paper have a number of practical implications for working with Archimedean copulas. We confine ourselves to brief comments on implications for copula construction, stochastic simulation and statistical inference.

5.1. *Construction of new distributions.* Many families of bivariate Archimedean copulas have been proposed, but models that extend to dimensions three and higher are much less common. The majority of such copula families have generators which are Laplace transforms of nonnegative random variables and, as discussed in Section 4.3, can only capture positive dependence. Examples of multivariate Archimedean copulas whose generators are not completely monotone seem to be rare; the standard examples are the Frank or the Clayton family introduced in Example 2.3. Theorem 3.1 offers an elegant solution to this problem: *Any $d$-dimensional Archimedean copula can be obtained simply as a survival copula of a $d$-dimensional $\ell_1$-norm* symmetric distribution with some radial part $R$. Concretely, by selecting a $p$-parametric family $\mathcal{F}_\theta = \{F_\theta : \theta \in \Theta \subseteq \mathbb{R}^p\}$ of distributions on $[0, \infty)$, we obtain a parametric family $\mathcal{C}_\theta$ of $d$-dimensional Archimedean copulas whose generators are given by

$$\{\psi_\theta(x) = \mathfrak{W}_d F_\theta(x) : \theta \in \Theta\}.$$

EXAMPLE 5.1. Consider the two parametric family of distributions with densities

$$f_{a,b}(x) = \begin{cases} \dfrac{ab}{b-a} x^{-2}, & \text{if } a \leq x \leq b, \\ 0, & \text{otherwise,} \end{cases}$$

where $0 < a < b$. $f_{a,b}(x)$ is the density of the reciprocal of a uniform distribution on the interval $[b^{-1}, a^{-1}]$. The survival copula of the $\ell_1$-norm symmetric distribution with this radial distribution has generator

$$\begin{aligned} \psi_{a,b}(x) &= \mathfrak{W}_d F_{a,b}(x) \\ &= \int_a^b \frac{ab}{b-a} t^{-2} \left(1 - \frac{x}{t}\right)_+^{d-1} dt \\ &= \frac{ab}{xd(b-a)} \left(\left(1 - \frac{x}{b}\right)_+^d - \left(1 - \frac{x}{a}\right)_+^d\right). \end{aligned}$$

Clearly, $\psi_{a,b}(x)$ is not completely monotone because it does not have derivatives of all orders. In Figure 4 we show four examples of 2000 points generated from the corresponding bivariate copula for $a = 1$ and different values of $b$.

5.2. *Stochastic simulation.* For a given Archimedean generator $\psi \in \Psi_d$ generating a copula $C$ in dimension $d$, it is of interest to have a general algorithm for generating random vectors $\mathbf{U}$ with distribution function $C$. If $\psi \in \Psi_\infty$, then a known method of generation is to use the mixed exponential or frailty representation discussed at the beginning of Section 3. This idea can be traced back to [23]; see also [24]. A key feature of the method is that we need to invert the Laplace transform $\psi$ to find a probability distribution function $F_W$ on the positive half-axis, from which we must then be able to sample; this is not always straightforward.

The insights of Section 3 give an alternative method. The steps of this method may also be found in [39], but are only justified there for generators $\psi \in \Psi_\infty$. In fact, the method may in principle be applied to any generator $\psi \in \Psi_d$ for $d \geq 2$. To generate a random vector from the $d$-dimensional copula $C$ with generator $\psi$, the algorithm is as follows:

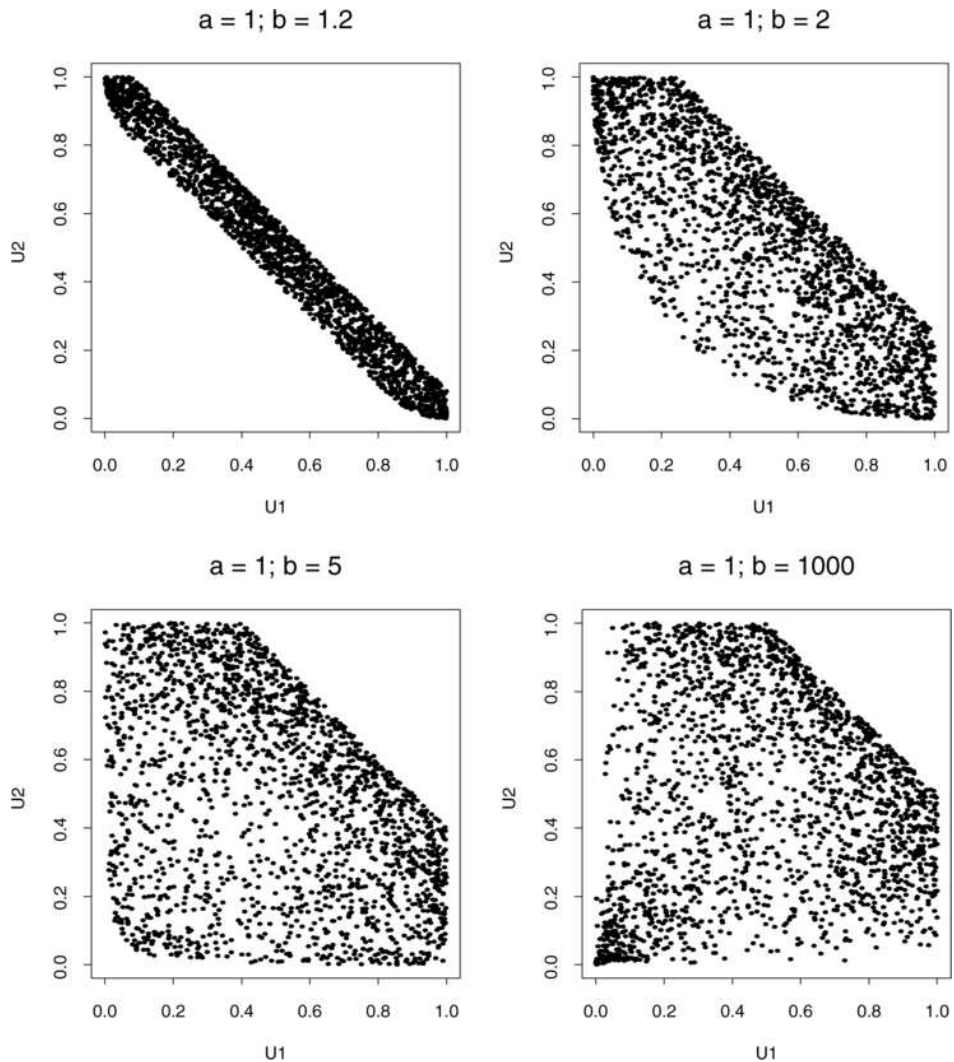

Fig. 4. *Four examples of 2000 points generated from the copula of Example 5.1 for $a = 1$ and different values of $b$.*

1. Generate a random vector $\mathbf{S} = (S_1, \ldots, S_d)$ uniformly distributed on the $d$-dimensional simplex $\mathcal{S}_d$. Based on Example 3.2, we can generate a vector of i.i.d. standard exponential variates $\mathbf{Y} = (Y_1, \ldots, Y_d)$ and return $S_i = Y_i / \|\mathbf{Y}\|_1$ for $i = 1, \ldots, d$.
2. Generate a univariate random variable $R$ having the radial distribution $F_R$ associated with $\psi$ in the sense of Definition 3.3; this will generally require us to compute the inverse Williamson $d$-transform $F_R = \mathfrak{W}_d^{-1} \psi$.
3. Return $\mathbf{U} = (U_1, \ldots, U_d)$ where $U_i = \psi(RS_i)$ for $i = 1, \ldots, d$.

The second step is the practically challenging part of the algorithm. In general we would use the inversion method to generate $R$ although this may require both accurate evaluation of $F_R = \mathfrak{W}_d^{-1}\psi$ and numerical inversion of $F_R$; see Example 3.3 for the distribution function that would have to be computed and inverted for the Clayton copula. An example where the method works particularly efficiently and allows us to generate samples from an attractive bivariate family is given below.

EXAMPLE 5.2. Consider the generator $\psi(t) = (1 - t^{1/\theta})_+$ where $\theta \geq 1$; this is the second family considered in [27] (see Table 4.1). Since $\psi$ is convex, it follows that $\psi \in \Psi_2$ but $\psi \notin \Psi_3$ because $\psi^{(1)}(t)$ is not defined at $t = 1$. Since $\psi_+^{(1)}(t) = -\theta^{-1} t^{1/\theta - 1}$ for $t \in [0, 1)$ and $\psi_+^{(1)}(t) = 0$ otherwise, we can use Theorem 3.1, part (ii), to show that

$$F_R(x) = \begin{cases} \left(1 - \dfrac{1}{\theta}\right) x^{1/\theta}, & \text{if } 0 \leq x < 1, \\ 1, & \text{if } x \geq 1, \end{cases}$$

which has a jump at $x = 1$ and yields a bivariate copula with a singular component. The distribution of $R$ is easily sampled by inversion. In Figure 5 we show four examples of 2000 points generated from this copula for different values of $\theta$. Obviously, $\theta = 1$ corresponds to the Fréchet–Hoeffding lower bound copula generated by $\psi_2^{\mathbf{L}}$ and, as $\theta \to \infty$, we approach perfect positive dependence. It is easy to use Proposition 4.7 to calculate that the Kendall's rank correlation of this copula is given by $\tau(C) = 1 - 2/\theta$, which implies that the points in the second picture ($\theta = 2$) are taken from a copula with a Kendall's rank correlation equal to zero.

5.3. *Statistical inference.* Proposition 4.4 can be used to devise diagnostic tests for particular hypothesized copulas. Suppose we have $d$-dimensional data vectors $\mathbf{U}_1, \ldots, \mathbf{U}_n$ that are believed to form a random sample from an Archimedean copula with generator $\psi$. (Alternatively they might be observations from the empirical copula of some multivariate data, derived by applying the marginal empirical distribution functions to the data as in [10].) For each vector observation $\mathbf{U}_k = (U_{k1}, \ldots, U_{kd})$ we can form the quantities $R_k = \sum_{i=1}^d \psi^{-1}(U_{ki})$ and $\mathbf{S}_k = (\psi^{-1}(U_{k1}), \ldots, \psi^{-1}(U_{kd}))/R_k$.

The $R_k$ data should be compatible with the radial distribution $F_R = \mathfrak{W}_d^{-1}\psi$, whereas the $\mathbf{S}_k$ data should be compatible with the hypotheses of uniformity on $\mathcal{S}_d$ and independence from the radial parts. Clearly, there are various possible ways of deriving numerical or graphical tests of these hypotheses.

It is also possible to reduce the estimation of Archimedean copulas to a one-dimensional problem by focusing on radial distributions. Ideas developed by [15] for the bivariate case may be extended to the general multivariate case. In that paper a nonparametric margin-free estimation procedure is proposed for the distribution function $K_C(x) = P(C(U_1, U_2) \leq x)$, where $(U_1, U_2)$ are random variables distributed according to the unknown copula $C$. This procedure may be applied in the higher-dimensional case to derive a nonparametric estimate of $K_C(x) = P(C(U_1, \ldots, U_d) \leq x)$, which is also given by $K_C(x) = P(R \geq \psi^{-1}(x))$, where $R$ is the radial part of the $\ell_1$-norm symmetric distribution associated with the generator $\psi$ of $C$. Thus the distribution $K_C$ may also be computed under parametric assumptions for the copula generator $\psi$ using Proposition 4.5. This makes it possible to estimate the copula by choosing the parameter values that give the best correspondence between the parametric model for $K_C$ and the nonparametric estimate.

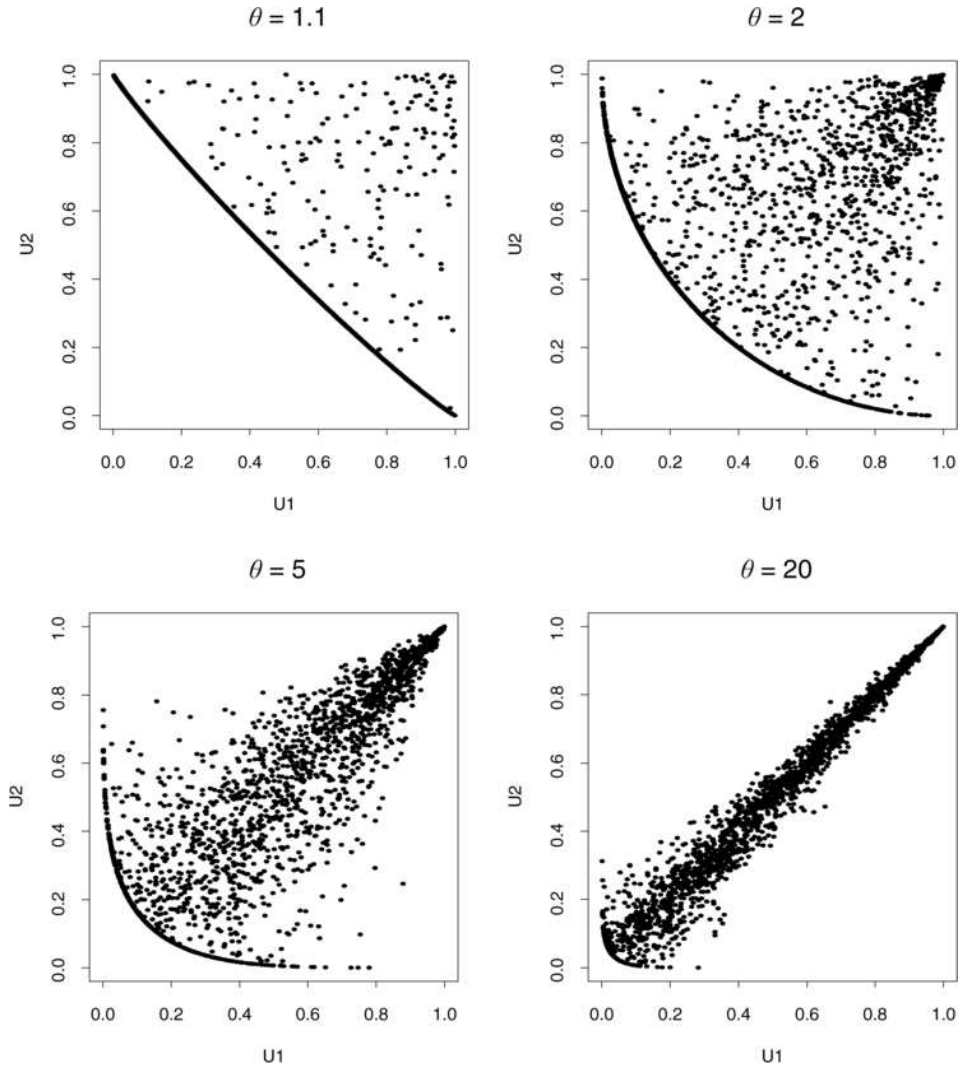

Fig. 5. *Four examples of 2000 points generated from the copula of Example 5.2 for different values of $\theta$.*

## APPENDIX A: PROOFS FROM SECTION 2

Theorem 2.2 relies on several results on $d$-monotone functions. First is the following lemma which summarizes Theorems 4, 5 and 6 of Chapter IV of [40].

LEMMA 2. *Let $d \geq 1$ be an integer and $f$ a nonnegative real function on $(a, b)$, $a, b \in \overline{\mathbb{R}}$. If $f$ satisfies*

(A.1) $$(\Delta_h)^k f(x) \geq 0$$

*for any $k = 1, \ldots, d$, any $x \in (a, b)$ and any $h > 0$ so that $x + kh \in (a, b)$, then*

(i) *$f$ is nondecreasing on $(a, b)$.*
(ii) *If $d \geq 2$, $f$ is convex and continuous on $(a, b)$.*
(iii) *If $d \geq 2$, the left-hand and right-hand derivatives of $f$ exist everywhere in $(a, b)$. Moreover, for any $a < x < y < b$,*

$$f'_-(x) \leq f'_+(x) \leq f'_-(y).$$

(iv) *For any $k = 1, \ldots, d-1$ and any $a < x \leq y < b$, $(\Delta_h)^k f(x) \leq (\Delta_h)^k f(y)$ whenever $h$ is chosen small enough so that $y + kh < b$.*

The second result which has proved of key importance is the following characterization of $d$-monotone functions:

PROPOSITION A.1. *Let $f$ be a real function on $(a, b)$, $a, b \in \overline{\mathbb{R}}$ and let $\tilde{f}$ denote a function on $(-b, -a)$ given by $\tilde{f}(x) = f(-x)$. Further, let $d \geq 1$ be an integer. Then the following statements are equivalent:*

(i) *$f$ is $d$-monotone in $(a, b)$.*
(ii) *$f$ is nonnegative and satisfies, for any $k = 1, \ldots, d$, any $x \in (-b, -a)$ and any $h_i > 0$, $i = 1, \ldots, k$ so that $(x + h_1 + \cdots + h_k) \in (-b, -a)$,*

$$\Delta_{h_k} \cdots \Delta_{h_1} \tilde{f}(x) \geq 0,$$

*where $\Delta_{h_k} \cdots \Delta_{h_1}$ denotes a sequential application of the first-order difference operator $\Delta_h$ given by $\Delta_h \tilde{f}(x) \equiv \tilde{f}(x + h) - \tilde{f}(x)$ whenever $x, x + h \in (-b, -a)$.*
(iii) *$f$ is nonnegative and satisfies, for any $k = 1, \ldots, d$, any $x \in (-b, -a)$ and any $h > 0$ so that $x + kh \in (-b, -a)$,*

$$(\Delta_h)^k \tilde{f}(x) \geq 0,$$

*where $(\Delta_h)^k$ denotes the $k$-monotone sequential use of the operator $\Delta_h$.*

Proof.

Proof of (ii) ⇔ (iii). First note that (iii) is a special case of (ii). For the reverse implication, assume without the loss of generality that $d \geq 2$ and fix an arbitrary $x$ in $(-b, -a)$ and $h_i > 0$, $i = 1, \ldots, d$ such that $x + h_1 + \cdots + h_d \in (-b, -a)$. Now, let $\tilde{f}_\ell$, $\ell = 1, \ldots, d$, denote a function on $(-b, -a - (h_d + \cdots + h_\ell))$ given by

$$\tilde{f}_\ell(y) = \Delta_{h_d} \cdots \Delta_{h_\ell} \tilde{f}(y), \qquad y \in (-b, -a - (h_d + \cdots + h_\ell)).$$

Observe in particular that any $\tilde{f}_\ell$, $\ell = 1, \ldots, d$, is nonnegative. One can even show that the functions $\tilde{f}_\ell$, $\ell = 2, \ldots, d$, satisfy (A.1) of Lemma 2 for any $k = 1, \ldots, \ell - 1$.

This claim can be verified by induction: because $\tilde{f}$ satisfies the assumptions of Lemma 2, the assertion (iv) thereof yields that $(\Delta_h)^k \tilde{f}(y) \leq (\Delta_h)^k \tilde{f}(y + h_d)$ for any $k = 1, \ldots, d - 1$ whenever $h$ is small enough so that $y + kh \in (-b, -a - h_d)$. But

$$(\Delta_h)^k \tilde{f}(y + h_d) - (\Delta_h)^k \tilde{f}(y) = (\Delta_h)^k \tilde{f}_d(y)$$

for any $k = 1, \ldots, d - 1$ and hence the beginning of the induction is established.

Now suppose that $\tilde{f}_{\ell+1}$ satisfies (A.1) for any $k = 1, \ldots, \ell$. The fact that $\tilde{f}_\ell(y) = \tilde{f}_{\ell+1}(y + h_\ell) - \tilde{f}_{\ell+1}(y)$ for any $y \in (-b, -a - (h_d + \cdots + h_\ell))$ yields

$$(\Delta_h)^k \tilde{f}_\ell(y) = (\Delta_h)^k \tilde{f}_{\ell+1}(y + h_\ell) - (\Delta_h)^k \tilde{f}_{\ell+1}(y)$$

for any $y \in (-b, -a - (h_d + \cdots + h_\ell))$ and any $k = 1, \ldots, \ell - 1$. The right-hand side is, however, nonnegative by (iv) of Lemma 2. This immediately yields that $\tilde{f}_\ell$ satisfies (A.1) for any $k = 1, \ldots, \ell - 1$.

Now, application of Lemma 2 to $\tilde{f}_2$ implies that $\tilde{f}_2$ is nondecreasing on $(-b, -a - (h_d + \cdots + h_2))$. Because $x, x + h_1 \in (-b, -a - h_d - \cdots - h_2)$ by assumption and

$$\tilde{f}_2(x + h_1) - \tilde{f}_2(x) = \Delta_{h_d} \cdots \Delta_{h_2} \tilde{f}(x + h_1) - \Delta_{h_d} \cdots \Delta_{h_2} \tilde{f}(x)$$
$$= \Delta_{h_d} \cdots \Delta_{h_1} \tilde{f}(x),$$

it follows that $\Delta_{h_d} \cdots \Delta_{h_1} \tilde{f}(x) \geq 0$. Since $\Delta_{h_k} \cdots \Delta_{h_1} \tilde{f}(x) \geq 0$ can be shown along the same lines for any $k = 1, \ldots, d - 1$, the desired implication follows. □

Proof of (i) ⇔ (iii). First recall that if $d = 1$ and $d = 2$, respectively, $d$-monotone monotonicity of $f$ reduces to $f$ nonincreasing and nonincreasing and convex, respectively. It is therefore sufficient to restrict the discussion to the case $d > 2$.

That (i) implies (ii) can be established by the midpoint theorem. Note first that if $f$ is $d$-monotone on $(a, b)$, then, for $k = 1, \ldots, d - 2$, $\tilde{f}^{(k)}$ exists on $(-b,-a)$ and is nonnegative there. Consequently, for any $k=1,\ldots,d-2$ and any $x \in (-b-a)$ and any $h > 0$ so that $x+kh \in (-b,-a)$ there exists $x^* \in (x, x+kh)$ so that $(\Delta_h)^k \tilde{f}(x) = \tilde{f}^{(k)}(x^*) \geq 0$.

The reverse implication can be established using the same argument as Theorem 7, Chapter IV of [40]. To do so, pick an arbitrary $x \in (-b,-a)$. The assertion (iv) of Lemma 2 yields that, for any $k=1,\ldots,d-1$ and any $h > 0$ such that $x+kh \in (-b,-a)$, $(\Delta_h)^k(\tilde{f}(x) - \tilde{f}(x-\ell)/\ell) \geq 0$ whenever $\ell > 0$ is sufficiently small so that $x - \ell \in (-b,-a)$. By letting $\ell \to 0$, it then follows that the left-hand derivative $\tilde{f}'_-$ fulfills (A.1) for any $k=1,\ldots,d-1$. The assertion (ii) of Lemma 2 in particular gives that $\tilde{f}'_-$ is continuous on $(-b,-a)$. Yet another application of Lemma 2, this time of the statement (iii) on $\tilde{f}$, then guarantees the existence of $\tilde{f}'$ on $(-b,-a)$.

The same chain of arguments can be applied successively on $\tilde{f}^{(k)}$. In the last step, this yields that $\tilde{f}^{(d-2)}$ exists on $(-b,-a)$ and fulfills (A.1) for $k=1,2$. In particular, $\tilde{f}^{(d-2)}$ is continuous, nondecreasing and convex. Put together, the derivatives $\tilde{f}^{(k)}$, $k=1,\ldots,d-2$, exist on $(-b,-a)$ and are continuous, nonnegative, nondecreasing and convex there. Because $(-1)^k f^{(k)}(x) = \tilde{f}^{(k)}(-x)$, the proof is complete. $\square$

REMARK A.1. After finishing this manuscript, Alfred Müller (personal communication) brought the notion of $d$-convexity as defined and studied by [18, 31, 32] to our attention. It is not our aim to discuss the various concepts of higher order convexity here; nonetheless, we would like to mention the paper [3] containing in particular the following results: First, a continuous function $f$ on $(a,b)$ which is *weakly $d$-convex* [meaning that $(\Delta_h)^d f(x) \geq 0$ for any $h > 0$ so that $x+dh \in (a,b)$] fulfills $\Delta_{h_1} \cdots \Delta_{h_d} f(x) \geq 0$ whenever $x + \sum_{i=1}^d h_i \in (a,b)$. Second, if a function $f$ is continuous and weakly $d$-convex on $(a,b)$, $f^{(d-2)}$ exists and is convex on $(a,b)$. These results would offer an alternative proof to parts of Proposition A.1. However, we favor our proof of the latter proposition, which is self-contained and more accessible. Furthermore, Proposition A.1 relates $d$-monotonicity, which is a stronger concept than weak $d$-convexity, to nonnegativity of higher order differences. This link is absolutely central to the study of quasi-monotonicity of survival functions as delineated in the proofs of Proposition 2.2 and Theorem 2.2 below.

The fact that Proposition A.1 characterizes $d$-monotonicity in terms of $\tilde{f}$ rather than $f$ may seem somewhat less elegant. The reason for this is that $d$-monotonicity is related to the existence of survival rather than distribution functions as described in Proposition 2.2, which we prove next.

PROOF OF PROPOSITION 2.2. Assume first that there exists a random vector $\mathbf{X}$ on $\mathbb{R}^d$ so that $\bar{H}(\mathbf{x}) = P(\mathbf{X} > \mathbf{x})$ for any $\mathbf{x} \in \mathbb{R}^d_+$. It is clear that $\lim_{x \to \infty} f(x) = 0$ and $f(0) = f(0+) = P(\mathbf{X} > \mathbf{0}) \in [0,1]$. Now, define $\tilde{f}(x)$ by $f(-x)$ for any $x \in (-\infty, 0]$ as in Proposition A.1 and observe that for any $x \in (-\infty, 0)$, any $k = 1, \ldots, d$ and any $h > 0$ so that $x + kh < 0$,

$$(\Delta_h)^k \tilde{f}(x) = \sum_{\ell=0}^{k} (-1)^{k-\ell} \binom{k}{\ell} \tilde{f}(x + \ell h)$$

$$= \sum_{\ell=0}^{k} (-1)^{k-\ell} \sum_{J \subseteq \{1,\ldots,k\}, \#J = \ell} f\left(-(k-\ell)\frac{x}{k} - \ell\left(\frac{x}{k} + h\right)\right)$$

$$= \sum_{\ell=0}^{k} (-1)^{k-\ell} \sum_{J \subseteq \{1,\ldots,k\}, \#J = \ell} P\left(\bigcap_{i \in J}\left\{-X_i < \frac{x}{k} + h\right\} \bigcap_{i \notin J}\left\{-X_i < \frac{x}{k}\right\}\right)$$

$$= P\left(\bigcap_{i=1}^{k}\left\{X_i \in \left(-\frac{x}{k} - h, -\frac{x}{k}\right]\right\}\right).$$

The last expression is, however, clearly nonnegative, in which case Proposition A.1 implies that $f$ is $d$-monotone on $(0, \infty)$. Because $f$ is right-continuous at 0, it is $d$-monotone even on $[0, \infty)$.

To establish the converse, note that $f$ is in particular continuous on $[0, \infty)$ by Lemma 2. Therefore, $\bar{H}$ is right-continuous. Because condition (i) of Lemma 1 is clearly fulfilled as well, the only assertion which needs to be verified is the quasi-monotonicity of $\bar{H}(-\mathbf{x})$. Since $\bar{H}$ is right-continuous and places no mass outside of $[0, \infty)^d$ and because $\mathbf{1}\{\mathbf{x} < \mathbf{0}\}$ is quasi-monotone, this in turn amounts to the confirmation of $\Delta_{\mathbf{h}} f(-\|\mathbf{x}\|_1) \geq 0$ for any $\mathbf{x} \in (-\infty, 0)^d$ and $\mathbf{h} > \mathbf{0}$ so that $\mathbf{x} + \mathbf{h} < \mathbf{0}$. However, this is immediately guaranteed by (ii) of Proposition A.1 upon setting $x = \|\mathbf{x}\|_1$. Hence, the proof is complete. □

With Proposition 2.2, the proof of Theorem 2.2 is straightforward.

PROOF OF THEOREM 2.2. As already discussed in Section 3, $C$ always satisfies the boundary conditions (i) and (ii) of Definition 2.1. Therefore, what remains to be shown is that $C$ is quasi-monotone if and only if $\psi$ is $d$-monotone on $(0, \infty)$. However, because any Archimedean generator is continuous and fulfills $\psi(0) = 1$, $\bar{F}$ given by $\bar{F}(x) = \psi(x)$ for $x \geq 0$ and by $\bar{F}(x) = 1$ otherwise, is a continuous univariate survival function.

Suppose first $\psi$ is $d$-monotone on $(0, \infty)$. Then, by Proposition 2.2, the function $\bar{H}$ defined as

$$\bar{H}(\mathbf{x}) = \psi(\|\max(\mathbf{x}, \mathbf{0})\|_1)$$

for any $\mathbf{x} \in \mathbb{R}^d$ is a survival function in $\mathbb{R}^d$ with continuous (identical) margins $\bar{F}$. Theorem 2.1 then ensures that the function $C$ associated with $\bar{H}$

via

$$C(\mathbf{u}) = \bar{H}(\bar{F}^{-1}(u_1), \ldots, \bar{F}^{-1}(u_d))$$

for any $\mathbf{u} \in [0,1]^d$ is a copula. Furthermore, $\bar{F}^{-1}(u_i) = \inf\{x : \psi(x) \leq u\} = \psi^{-1}(u_i)$ for $u_i \in [0,1)$ and $\bar{F}^{-1}(1) = \inf\{x : F(x) \geq 0\} = -\infty$. However, $\max(\bar{F}^{-1}(1), 0) = 0 = \psi^{-1}(1)$ and hence

$$C(\mathbf{u}) = \psi(\psi^{-1}(u_1) + \cdots + \psi^{-1}(u_d)), \qquad \mathbf{u} \in [0,1]^d.$$

Conversely, assume that $C$ is a copula. By virtue of Theorem 2.1 it therefore follows that $\bar{H}(\mathbf{x}) = C(\bar{F}(x_1), \ldots, \bar{F}(x_d))$, $\mathbf{x} \in \mathbb{R}^d$, is a joint survival function with (identical) marginal survival functions $\bar{F}$. Furthermore, it is easy to convince oneself that $\bar{H}(\mathbf{x}) = \psi(x_1 + \cdots + x_d)$ for any $\mathbf{x} \in \mathbb{R}^d_+$. Therefore, $\bar{H}$ satisfies the hypothesis of Proposition 2.2 and it follows that $\psi$ is $d$-monotone on $[0,\infty)$. □

## APPENDIX B: PROOFS FROM SECTION 3

PROOF OF PROPOSITION 3.1. Throughout, $(\mathscr{S})$ indicates Stieltjes integration; any other integral is to be understood in the usual Lebesgue–Stieltjes sense. The original result by [41] states that a function $f$ is $d$-monotone on $(0,\infty)$ if and only if there exists a nonnegative, nondecreasing real function $\gamma$ on $[0,\infty)$ which is bounded below and such that

$$f(t) = (\mathscr{S}) \int_0^\infty (1 - ut)_+^{d-1} \, d\gamma(u)$$

for any $t \in (0,\infty)$. Furthermore, $\gamma$ is uniquely determined at its points of continuity and fulfills $\gamma(0) = 0$ as well as $\gamma(0+) = f(\infty)$.

PROOF OF (i). Suppose $f$ is a Williamson $d$-transform of a nonnegative random variable $X$ with distribution function $F$. Then $f(0+) = f(0) = 1 - F(0)$ and $\lim_{x \to \infty} f(x) = 0$ by dominated convergence. That $f$ is $d$-monotone on $(0,\infty)$ is then clear from Williamson's result.

To establish the reverse implication, observe first that in the situation of Proposition 3.1, Williamson's result ensures the existence of a unique right-continuous, nonnegative, nondecreasing function $G$ which is bounded below and on $\mathbb{R}$ so that

$$f(t) = \int_0^\infty (1 - ut)_+^{d-1} \, dG(u)$$

for any $t \in (0,\infty)$ and $G(0-) = G(0+) = 0$ as well as $\lim_{u \to \infty} G(u) = p$. The transformation formula for Lebesgue–Stieltjes integrals then yields

$$f(t) = \int_0^\infty \left(1 - \frac{t}{u}\right)_+^{d-1} d\nu(u),$$

where $\nu$ is an image measure of the measure induced by $G$ with respect to the mapping $\phi$ given by $\phi(u) = 1/u$ for $u \in (0, \infty)$ and by 0 otherwise. A possible distribution function of $\nu$ is

$$\tilde{G}(u) = \begin{cases} 0, & \text{if } u \in (-\infty, 0), \\ p - G(1/u-), & \text{if } u \in [0, \infty). \end{cases}$$

The disturbing fact that $\nu$ is not necessarily a probability measure can be easily remedied upon setting

$$F(u) = \begin{cases} 0, & \text{if } u \in (-\infty, 0), \\ 1 - G(1/u-), & \text{if } u \in [0, \infty). \end{cases}$$

It then still holds that $f(t) = \int_0^\infty (1 - \frac{t}{u})_+^{d-1} dF(u)$ whenever $t > 0$ but also that $f(t) = \mathfrak{W}_d F(t)$ for any $t \in (0, \infty)$. Since both $f$ and $\mathfrak{W}_d F$ are right-continuous at 0, the reversed implication is established. $\square$

PROOF OF (ii). The fact that the relationship between $f$ and $F$ is one-to-one is clear from the original result by Williamson in combination with (i). Theorem 3 of [41], which is the inversion formula for Williamson $d$-transforms, further ensures that $F$ is of the form (3.3) at all points of continuity. Furthermore, the one-sided derivatives of $(-1)^{d-2} f^{(d-2)}$ exist everywhere and are nondecreasing by Lemma 2. A version of the midpoint theorem for continuous functions with one-sided derivatives (see [25]) easily yields that $f_+^{(d-1)}$ is right-continuous on $(0, \infty)$ with left-hand limit $f_-^{(d-1)}$. $\square$

## APPENDIX C: PROOFS FROM SECTION 4

PROOF OF PROPOSITION 4.1. Throughout, let $\mathbf{X} \sim H$ and $\mathbf{U} \sim C$ be random vectors and $\lambda_d$ be the Lebesgue measure on $\mathbb{R}^d$.

PROOF OF (i). It is sufficient to consider $d = 2$. In the case $d \geq 3$ the situation is simpler because then $\psi'$ exists everywhere in $(0, \infty)$ and one can argue by virtually the same arguments. Now, suppose $H$ has a density. Then obviously $P(\mathbf{X} \in (0, \psi^{-1}(0))^d) = 1$. Furthermore, consider the transformation $T : (0, \psi^{-1}(0))^d \to \mathbb{R}^d$ given by

$$T(x_1, \ldots, x_d) = (\psi(x_1), \ldots, \psi(x_d)), \qquad \mathbf{x} \in (0, \psi^{-1}(0))^d.$$

Then $T$ is injective and $\mathbf{U} \stackrel{\mathrm{d}}{=} T(\mathbf{X})$. Because $\psi'$ exists a.e. in $(0, \psi^{-1}(0))$, $T$ is a.e. regular (meaning that the set of points where $T$ is not differentiable has Lebesgue measure zero). Furthermore, if $\psi'(t)$ exists, then necessarily $\psi' > 0$.

To see this, assume for the moment the converse. Because $\psi$ is convex and decreasing on $(0, \psi^{-1}(0))$, $\psi'_+$ is nonpositive and nondecreasing. If now $0 < t < \psi^{-1}(0)$ were such that $\psi'(t) = 0$, there would exist an $\varepsilon > 0$ so that $\psi'(x) = 0$ for almost all $x \in [t, \psi^{-1}(0) - \varepsilon]$. However, $\psi$ is absolutely continuous by assumption so $\psi'(x) = 0$ a.e. in $[t, \psi^{-1}(0) - \varepsilon]$ would imply that $\psi$ is constant there, which is obviously not the case.

Put together, $D = \{\mathbf{u} \in (0,1)^d : \prod_{i=1}^d |\psi'(\psi^{-1}(u_i))| = 0\}$ is an empty set. In particular, $\lambda_d(D) = 0$ meaning that $T$ is a.e. smooth. By means of Sard's theorem and the second transformation theorem (see [17], Sections 8.9 and 8.10), it then follows that $C$ has a density.

The reverse statement can be established similarly by considering a transformation $\tilde{T}(0,1)^d \to \mathbb{R}^d$ given by

$$\tilde{T}(u_1, \ldots, u_d) = (\psi^{-1}(u_1), \ldots, \psi^{-1}(u_d)), \qquad \mathbf{x} \in (0,1)^d.$$

Clearly, $P(\mathbf{U} \in (0,1)^d) = 1$ and $\mathbf{X} \stackrel{d}{=} \tilde{T}(\mathbf{U})$ (for the latter claim, see the proof of Theorem 2.2). Again, one can easily convince oneself that $\tilde{T}$ is a.e. regular. Furthermore, $\psi^{-1}$ is (strictly) increasing, continuous and convex on $(0,1)$ and therefore in particular absolutely continuous in any $[a,b] \subset (0,1)$. By the same argument as above, it can then be established that $\psi^{-1'}(t) > 0$ whenever it exists. Consequently, $\tilde{T}$ is a.e. smooth and $\mathbf{X}$ has a density again by Sard's theorem and the second transformation theorem. $\square$

PROOF OF (ii). Suppose now $\psi \in \Psi_{d+1}$ and let $\tilde{C}$ be the $(d+1)$-dimensional Archimedean copula with generator $\psi$. Furthermore, fix a random vector $\tilde{\mathbf{X}}$ following the $\ell_1$-norm symmetric distribution associated with $\tilde{C}$ in the sense of Theorem 3.1. Arguments detailed in Section 5.2.3 in [8] then yield that all lower-dimensional margins of $\tilde{\mathbf{X}}$ have densities. In particular, this applies to $\mathbf{Y} = (\tilde{X}_1, \ldots, \tilde{X}_d)$. However, $\mathbf{Y} \stackrel{d}{=} \mathbf{X}$ and $C$ has a density by means of (i). $\square$

PROOF OF PROPOSITION 4.2. Because of (i) of Proposition 4.1 and (iii) of Lemma 2, we only need to verify that the radial part of the $\ell_1$-norm symmetric distribution associated with $C$ has a density if and only if $\psi_+^{(d-1)}$ is absolutely continuous on $(0, \infty)$. Before we do so, observe that for $x \in [0, \infty)$, $F_R$ may be rewritten as

$$F_R(x) = F_Q(x) - \frac{(-1)^{d-1} x^{d-1} f_+^{(d-1)}(x)}{(d-1)!},$$

where $F_Q$ refers to the distribution function of the radial part of the $\ell_1$-norm distribution associated with the $(d-1)$-dimensional margin of $C$, that is, with $\psi(\psi^{-1}(u_1) + \cdots + \psi^{-1}(u_{d-1}))$. Furthermore, Proposition 4.1 guarantees that $F_Q$ is absolutely continuous.

Assume first that $F_R$ is absolutely continuous and consider $[a, b] \subset (0, \infty)$. Because $f_+^{(d-1)}(x) = (-1)^{d-1}(d-1)! \frac{(F_Q(x) - F_R(x))}{x^{d-1}}$ on $[a, b]$, $f_+^{d-1}$ is absolutely continuous there.

To establish the reverse implication, observe that because $f_+^{d-1}$ is continuous on $(0, \infty)$, $F_Q$ is continuous on $\mathbb{R}$ and $F_R$ is continuous at 0, then $F_R$ is continuous in $\mathbb{R}$. Now, fix $[a, b] \subset (0, \infty)$. By a similar argument as above it then clearly holds that $F_R$ is absolutely continuous in $[a, b]$. Because $F_R$ is continuous and of bounded variation, it is absolutely continuous even in $[0, b]$. Since $F_R$ is clearly absolutely continuous on $(-\infty, 0]$, the absolute continuity of $F_R$ on $\mathbb{R}$ is immediate.

The last claim follows from the fact that if $H$ is an absolutely continuous distribution function then its density $h$ satisfies

$$h(\mathbf{x}) = \frac{\partial^d}{\partial x_1 \cdots \partial x_d} H(\mathbf{x})$$

for almost all $\mathbf{x} \in \mathbb{R}^d$. However clear this statement may seem, it is far from obvious; see [34], page 115 or [6]. □

PROOF OF PROPOSITION 4.3. First recall that if $\mathbf{U} \sim C$, then the random vector $\mathbf{X} \stackrel{d}{=} (\psi^{-1}(U_1), \ldots, \psi^{-1}(U_d))$ follows an $\ell_1$-norm symmetric distribution associated with $C$. Denote the radial part of $\mathbf{X}$ by $R$ and its distribution function by $F_R$. Then, for any $s \in (0, 1]$,

$$P^C(L(s)) = P\left(\sum_{i=1}^d \psi^{-1}(U_i) = \psi^{-1}(s)\right)$$

$$= P\left(\sum_{i=1}^d X_i = \psi^{-1}(s)\right) = P(R = \psi^{-1}(s)),$$

where the last equality follows from (iii) of Proposition 3.2. Inversion (3.3) and the fact that $\psi^{(k)}$ is continuous for any $k = 1, \ldots, d-2$ in turn imply that

$$P(R = \psi^{-1}(s))$$
$$= F_R(\psi^{-1}(s)) - F_R(\psi^{-1}(s)-)$$
$$= \frac{(-1)^{d-1}(\psi^{-1}(s))^{d-1}}{(d-1)!}(\psi_-^{(d-1)}(\psi^{-1}(s)) - \psi_+^{(d-1)}(\psi^{-1}(s))).$$

Regarding $L(0)$, one similarly argues that $P^C(L(0)) = P(R \geq \psi^{-1}(0))$. From this it is first immediate that $P^C(L(0)) = 0$ whenever $\psi^{-1}(0) = \infty$. Suppose now $\psi^{-1}(0) < \infty$. Because all derivatives $\psi^{(k)}(x)$, $k = 1, \ldots, d-2$, as well as $\psi_+^{(d-1)}(x)$ vanish for $x \in [\psi^{-1}(0), \infty)$, it ensues from (3.3) that $F_R(x) = 1$ for any $x \in [\psi^{-1}(0), \infty)$. Therefore,

$$P^C(L(0)) = P\left(\sum_{i=1}^d \psi^{-1}(U_i) = \psi^{-1}(0)\right).$$

Along the same lines as in the case of $L(s)$, $s \in (0,1]$, one reasons that

$$P^C(L(0)) = P(R = \psi^{-1}(0)) = \frac{(-1)^{d-1}(\psi^{-1}(0))^{d-1}}{(d-1)!}\psi_-^{(d-1)}(\psi^{-1}(0)),$$

which concludes the proof. $\square$

PROOF OF PROPOSITION 4.6. First, recall that $C_d^{\mathbf{L}}$ is Archimedean, as detailed in Example 2.2. Furthermore, $\psi_2^{\mathbf{L}}(x) = (1-x)_+$ and $C_2^{\mathbf{L}}$ coincides with the Fréchet–Hoeffding lower bound. Therefore, assume that $d \geq 3$ and let $\psi$ be the generator of $C$. Because $(\psi_d^{\mathbf{L}})^{-1}(x) = (1 - x^{1/(d-1)})$, the assertion $C_d^{\mathbf{L}} \prec_{\mathrm{cL}} C$ can be shown by verifying that the function $f(x) = (\psi_d^{\mathbf{L}})^{-1}(\psi(x))$ is concave on $(0, \infty)$; see, for example, [27], Corollary 4.4.4. If $\psi^{-1}(0) < \infty$, it holds that $f(x) = 1$ for $x \in [\psi^{-1}(0), \infty)$ and it is sufficient to show that $f$ is concave on $[0, \psi^{-1}(0)]$.

Assume first $d \geq 4$. In that case, $f$ is twice differentiable on $(0, \psi^{-1}(0))$ with

$$f^{(2)}(x) = -\frac{1}{d-1}\psi^{(2)}(x)\psi(x)^{-(d-2)/(d-1)}$$

$$+ \frac{d-2}{(d-1)^2}(\psi'(x))^2\psi(x)^{-(d-2)/(d-1)-1}$$

$$= \frac{(d-2)(\psi'(x))^2 - (d-1)\psi^{(2)}(x)\psi(x)}{(d-1)^2\psi(x)^{(d-2)/(d-1)+1}}.$$

To establish that $f^{(2)}(x) \leq 0$ for $x \in (0, \psi^{-1}(0))$, we consider the matrix

$$A = \begin{pmatrix} \psi(x) & -\dfrac{1}{d-1}\psi'(x) \\ -\dfrac{1}{d-1}\psi'(x) & \dfrac{1}{(d-1)(d-2)}\psi^{(2)}(x) \end{pmatrix}$$

and show that it is positive semi-definite for all $x \in (0, \psi^{-1}(0))$. Since $\psi$ is $d$-monotone, the dominated convergence theorem ensures that, for $x \in (0, \psi^{-1}(0))$,

$$\psi^{(k)}(x) = (d-1)\cdots(d-k)\int_0^\infty (-1)^k \frac{1}{t^k}\left(1 - \frac{x}{t}\right)_+^{d-1-k} dF_R(t),$$

$$k = 1, \ldots, d-2,$$

where $F_R$ is the distribution function of a nonnegative random variable whose Williamson $d$-transform is $\psi$. Now, fix $a_i \in \mathbb{R}$, $i = 1, 2$. Then

$$\sum_{i=1}^{2}\sum_{i=1}^{2} a_i a_j A_{ij} = (a_1)^2 \int_0^\infty \left(1 - \frac{x}{t}\right)_+^{d-1} dF_R(t)$$

$$+ 2a_1 a_2 \int_0^\infty \frac{1}{t}\left(1 - \frac{x}{t}\right)_+^{d-2} dF_R(t)$$

$$+ (a_2)^2 \int_0^\infty \frac{1}{t^2}\left(1 - \frac{x}{t}\right)_+^{d-3} dF_R(t)$$

$$= \int_x^\infty \left(1 - \frac{x}{t}\right)^{d-3} \left(a_1\left(1 - \frac{x}{t}\right) + \frac{a_2}{t}\right)^2 dF_R(t) \geq 0.$$

Consequently, $|A| \geq 0$, which in turn implies $f^{(2)}(x) \leq 0$ for $x \in (0, \psi^{-1}(0))$.

For $d = 3$ we can use the fact that $\psi_+^{(2)}$ exists everywhere in $(0, \psi^{-1}(0))$. Consequently, $f_+^{(2)}(x)$ exists for any $x \in (0, \psi^{-1}(0))$ and

$$f_+^{(2)}(x) = \frac{(d-2)(\psi'(x))^2 - (d-1)\psi_+^{(2)}(x)\psi(x)}{(d-1)^2 \psi(x)^{(d-2)/(d-1)+1}}.$$

Concavity of $f$ then follows from $f_+^{(2)} \leq 0$ on $(0, \psi^{-1}(0))$; see Theorem 1 of [25]. Now consider $h > 0$ and $x \in (0, \psi^{-1}(0))$. Then

$$\frac{\psi'(x+h) - \psi'(x)}{h}$$
$$= -2\left[\frac{\int_{x+h}^\infty (1/t)(1-(x+h)/t)\,dF_R(t) - \int_x^\infty (1/t)(1-x/t)\,dF_R(t)}{h}\right]$$
$$= 2\left[\int_{(x+h,\infty)} \frac{1}{t^2}\,dF_R(t) + \int_{(x,x+h]} \frac{1}{th}\left(1 - \frac{x}{t}\right) dF_R(t)\right].$$

Because $\frac{1}{th}(1 - \frac{x}{t}) \leq \frac{1}{x(x+h)}$ for $t \in (x, x+h]$, the dominated convergence theorem yields

$$\psi_+^{(2)}(x) = \lim_{h \downarrow 0} \frac{\psi'(x+h) - \psi'(x)}{h} = 2 \int_{(x,\infty)} \frac{1}{t^2}\,dF_R(t).$$

The inequality $f_+^{(2)}(x) \leq 0$ for $x \in (0, \psi^{-1}(0))$ can now be established by exactly the same arguments as in the case $d \geq 4$. □

**Acknowledgments.** The authors would like to thank Professors Paul Embrechts, Christian Genest and Alfred Müller for fruitful discussions and an anonymous referee for useful comments.

DEPARTMENT OF ACTUARIAL MATHEMATICS
AND STATISTICS
HERIOT-WATT UNIVERSITY
EDINBURGH EH14 4AS
UNITED KINGDOM
E-MAIL: A.J.McNeil@hw.ac.uk
URL: http://www.ma.hw.ac.uk/~mcneil/

DEPARTMENT OF MATHEMATICS
ETH ZURICH
CH-8092 ZURICH
SWITZERLAND
E-MAIL: johanna@math.ethz.ch
URL: http://www.math.ethz.ch/~johanna/